\documentclass[11pt]{article}
\usepackage[numbers]{natbib}
\usepackage{amsbsy}
\usepackage{amsgen}

\usepackage{amstext}
\usepackage{amsxtra}
\usepackage{amsmath}
\usepackage{amsopn}
\usepackage{amsthm}
\newtheorem{thm}{Theorem}[section]
\theoremstyle{definition}
\newtheorem{dfn}{Definition}[section]
\theoremstyle{plain}

\theoremstyle{remark}
\newtheorem{note}{Remark}[section]
\theoremstyle{plain}
\newtheorem{lem}[thm]{Lemma}
\theoremstyle{plain}

\theoremstyle{plain}

\theoremstyle{conjecture}
\newtheorem{con}{Conjecture}[section]

\begin{document}\date{}

\title{SYZ Mirror Symmetry of Hitchin's Moduli Spaces Near Singular Fibers I}
\author{Wenxuan Lu\footnote{wenxuanl@math.upenn.edu}} \maketitle
\begin{abstract}We study hyperkahler metrics and hyperholomorphic connections of Hitchin's moduli spaces after Gaiotto, Moore and Neitzke. Their construction via the twistor technique produces intricate wall crossing behaviors. For certain four dimensional Hitchin's moduli spaces  local models and degeneration to local models near singular fibers of the Hitchin's fibration are understood. \end{abstract}
\newpage
\tableofcontents{\section{Introduction}

Hitchin's moduli spaces are very important in geometry, representation theory and mathematical physics. The Strominger-Yau-Zaslow type mirror symmetry of Hitchin's moduli spaces with Langlands dual gauge groups is believed to be able  to recover the geometric Langlands duality since the seminal work of \citep{HT} and \citep{KW}. On the other hand, in recent years Gaiotto, Moore and Neitzke have developed a formalism that can produce an indirect yet exact description of hyperkahler metrics and certain hyperholomorphic connections on $SU(2)$ Hitchin's moduli spaces in a series of papers \citep{GMN1,GMN2,GMN3}. Interesting questions emerge from these progresses. \begin{enumerate}
\item What does the SYZ mirror symmetry say about hyperkahler metrics and hyperholomorphic connections in the Gaiotto-Moore-Neitzke construction?  This is not only potentially important for the geometric Langlands program but also useful even for general studies of nontrivial\footnote{Here NONTRIVIAL means beyond the semiflat mirror symmetry.} SYZ mirror symmetry of Calabi-Yau spaces.
For example, how does the mirror duality act on the hyperkahler metric? Does it mean anything to the geometric Langlands duality? Another example is that we know mirror symmetry should interchange the so-called $A$ branes and $B$ branes. For Hitchin's moduli spaces we have a more sophisticated scenario: the correspondence of $(B,B,B)$ branes and $(B,A,A)$ branes\footnote{See \citep{KW} for the meaning of these terminologies.}. Hyperholomorphic connections are $(B,B,B)$ branes. What are the mirrors? What is the relation between these branes (and the mirror symmetry of them) and the geometric Langlands duality?
\item The general philosophy of SYZ mirror symmetry predicts that the construction of the mirror family of a given family of Calabi-Yau spaces should involve the so-called $instanton\ corrections$ from the counting of some geometric objects. The Gaiotto-Moore-Neitzke construction has the feature of instanton corrections and the wall crossing phenomenon is nothing but a consistency condition for contributions of these corrections. Further work \citep{L} shows the compatibility with the SYZ mirror symmetry picture. Then naturally one wants to ask: what is the role of wall crossing phenomenon in the geometric Langlands duality?
\item One of the major difficulties in the study of mirror symmetry is that we do not know how to handle singular fibers in a special Lagrangian fibration. For Hitchin's moduli spaces singular fibers of the so-called Hitchin's fibration should contain important information in representation theory. Frenkel and Witten have shown in \citep{FW} how to handle this singular fiber problem for the simplest branes and discovered that indeed it is related to a deep aspect of the  Langlands duality called endoscopy. We might wonder if we can deal with more complicated aspects of mirror symmetry near singular fibers and clarify their representation theoretical meanings.\end{enumerate}

This paper is a very modest step toward the study of these problems.  The basic idea is that if one has a good local model near certain singular fiber then one should be able to use it to approximate the full geometry of Hitchin's moduli spaces and extract some information from it. Hopefully the local information is simple enough that one can analyze (some local version of) mirror symmetry and has nontrivial meanings in the Langlands duality\footnote{Even local information can mean something quite nontrivial like the situation of Frenkel-Witten case.}.

In this paper we describe some nice local models called Ooguri-Vafa spaces and show that some four dimensional Hitchin's moduli spaces do degenerate to them near singular fibers. In particular we will see how the wall crossing behavior of the full geometry and the local model match. There is no essential obstacle of generalizing it to higher dimensional cases near singular fibers that can have the higher dimensional generalization of Ooguri-Vafa spaces as local models\footnote{However in higher dimensions this is probably not very useful because singular fibers can have much more complicated local models than Ooguri-Vafa spaces. The lack of more local models in differential geometry prevents us from revealing essentially new phenomena that can not be seen in the four dimensional case. That is why we will stick to the dimension 4.}. We also give a (partially conjectural) description of the metrical aspect of the local version of mirror symmetry. A more ambitious goal would be to study the SYZ mirror symmetry of branes. We plan to consider it in the future.

Section 2 is an exposition and clarification of Gaiotto, Moore and Neitzke's construction. We hope that this section and some sections of the author's previous paper \citep{L} can help mathematicians without physics background understand their work. Section 3 introduces the Hitchin's moduli spaces we want to study. They are the ones that Frenkel and Witten studied. Section 4 and 5 are materials about Ooguri-Vafa spaces. In section 6 we analyze the relation between local models and the full geometry of Hitchin's moduli spaces.\\

\noindent {\bf Acknowledgements} The author would like to thank Ron Donagi, Sean Keel and Andrew Neitzke for discussions on various related topics.
\section{Hyperkahler Metrics and Hyperholomorphic Connections on Hitchin's Moduli Spaces}
\subsection{Gaiotto-Moore-Neitzke 2d-4d Wall Crossing Formula}

The following formalism of wall crossing formulas is introduced by  Gaiotto-Moore-Neitzke in \citep{GMN1}.

Let $\mathcal{V}$ be a finite set and $\Gamma$ be a lattice with an integral antisymmetric pairing $\langle\cdot,\cdot\rangle$. Let $\Gamma_{i}$ be an $\Gamma$-torsor for each $i\in \mathcal{V}$. $\Gamma_{i}$ is assumed to be both a left and right $\Gamma$-torsor. We use $\gamma_{i}$ to denote an element of $\Gamma_{i}$. Define the set $\Gamma_{ij}$ as the set of formal differences $\gamma_{i}-\gamma_{j}$ modulo the relation $\gamma_{i}-\gamma_{j}=(\gamma_{i}+\gamma)-(\gamma_{j}+\gamma)$. $\Gamma_{ij}$ is also a left and right $\Gamma$-torsor.An element of $\Gamma_{ij}$ is denoted by $\gamma_{ij}$. $\Gamma_{ii}$ can be identified with $\Gamma$. There are natural addition operations $\Gamma_{ij}\times\Gamma_{j}\rightarrow \Gamma_{i}$ given by $\gamma_{ij}+\gamma_{j}$ and $\Gamma_{ij}\times\Gamma_{jk}\rightarrow \Gamma_{ik}$ given by $\gamma_{ij}+\gamma_{jk}$.

We assume the existence of a $central \ charge$ function $Z$. First we consider the lattice $\Gamma$.  $Z: \Gamma\rightarrow \mathbf{C}$ is a linear function on $\Gamma$. We use $Z_{\gamma}$ to denote $Z(\gamma)$. We then assume that it has an extension to all $\Gamma_{i}$ compatible with the $\Gamma$-action i.e. $Z_{\gamma+\gamma_{i}} = Z_{\gamma}+Z_{\gamma_{i}}$. This induces a further extension to $\Gamma_{ij}$ by defining $Z_{\gamma_{ij}}:=Z_{\gamma_{i}}-Z_{\gamma_{j}}$.

We assume the existence of a collection of numbers called $BPS \ numbers$. An integer $\Omega_{\gamma}$ for each $\gamma\in\Gamma$, an single rational number denoted by $\omega(\gamma,\gamma_{i})$ for each $i$ and an integer $\mu(\gamma_{ij})$ for each $\gamma_{ij}\in \Gamma_{ij}$ where $i\neq j$. Then we define $\omega: \Gamma\times\Gamma\rightarrow \mathbf{Z}$
\begin{equation}\omega(\gamma,\gamma^{'}):=\Omega(\gamma)\langle\gamma,\gamma^{'}\rangle\end{equation}We  define $\omega: \Gamma\times \coprod_{i}\Gamma_{i}\rightarrow \mathbf{Z}$ by \begin{equation}\omega(\gamma,\gamma_{i}+\gamma^{'}):=\omega(\gamma,\gamma_{i})+\Omega(\gamma)\langle\gamma,\gamma^{'}\rangle\end{equation}
We also define $\omega: \Gamma\times \coprod_{i,j}\Gamma_{ij}\rightarrow \mathbf{Z}$ by\begin{equation}\omega(\gamma,\gamma_{i}-\gamma_{j})=\omega(\gamma,\gamma_{ij}):=\omega(\gamma,\gamma_{i})-\omega(\gamma,\gamma_{j})\end{equation}
It is  required as a part of the definition of BPS numbers that $\omega(\gamma,\gamma_{ij})$ is always an integer.

We also assume the existence of a $twisting \ function$ $\sigma(a,b)$ valued in $\{\pm 1\}$where $a, b$ are in $\Gamma,\Gamma_{i}$ or $\Gamma_{ij}$. It is defined whenever $a+b$ is defined and satisfies \begin{equation}\sigma(a,b)\sigma(a+b,c)=\sigma(a,b+c)\sigma(b,c)\end{equation}We also assume that if $a+b$ and $b+a$ are both defined then $\sigma(a,b)=\sigma(b,a)$ and $\sigma(\gamma,\gamma^{'})=(-1)^{\langle\gamma,\gamma^{'}\rangle}$.

We introduce a set of formal variables $X_{a}$ where $a$ is in $\Gamma,\Gamma_{i}$ or $\Gamma_{ij}$. We define a multiplication \begin{equation}X_{a}X_{b}=\sigma(a,b)X_{a+b} \ if \ a+b\ is \ defined; X_{a}X_{b}=0\ otherwise\end{equation}

We define rays on the complex plane by \begin{equation}l_{\gamma}:=Z_{\gamma}\mathbf{R}_{-};\ \ \ l_{\gamma_{ij}}:=Z_{\gamma_{ij}}\mathbf{R}_{-}\end{equation} where $\mathbf{R}_{-}$ is the set of negative real numbers. $l_{\gamma}$ is called a $BPS \ \mathcal{K}-ray$ if $\omega(\gamma,\cdot)\neq 0$ and $l_{\gamma_{ij}}$ is called a $BPS\  \mathcal{S}-ray$ if $\mu(\gamma_{ij})\neq 0$.

For a BPS $\mathcal{S}$-ray we associate a transformation called the $\mathcal{S}-factor$ defined as \begin{equation}\mathcal{S}_{\gamma_{ij}}^{\mu}: X_{a}\rightarrow (1-\mu(\gamma_{ij})X_{\gamma_{ij}})X_{a}(1+\mu(\gamma_{ij})X_{\gamma_{ij}})\end{equation}

For a BPS $\mathcal{K}$-ray we associate a transformation called the $\mathcal{K}-factor$ defined as \begin{equation}\mathcal{K}_{\gamma}^{\omega}: X_{a}\rightarrow (1-X_{\gamma})^{-\omega(\gamma,\ a)}X_{a}\end{equation}

Let $\forall$ be a convex angular sector with the apex at the origin we define a transformation  \begin{equation}A(\forall):=\prod_{\gamma:\ l_{\gamma}\in\forall}\ \ \prod_{\gamma_{ij}:\ l_{\gamma_{ij}}\in\forall}\mathcal{K}_{\gamma}^{\omega}\mathcal{S}_{\gamma_{ij}}^{\mu}\end{equation}
The product here is understood as an ordered product. In other words (type $\mathcal{S}$ or $\mathcal{K}$) factors are ordered so that we encounter them following the order of the corresponding rays along the counterclockwise direction.

The order of those rays and BPS numbers depends on the choice of the central charge\footnote{We do not keep the dependence of BPS numbers on the central charge in our notations in this paper. But one should be aware of this fact (see \citep{GMN2},\citep{L} for further studies of this issue).}. The order does not  change if we take a small enough perturbation of the central charge. However if we vary the central charge $Z$ in the space of central charge functions and let it cross the real codimention one locus where two ($\mathcal{S}$ or $\mathcal{K}$) rays coincide in the complex plane then the order of rays in the definition of $A(\forall)$ will change. Such a real codimention one locus is called a $marginal\ stability\ wall$. We shall call it a $stability\ wall$ in this paper. Stability walls decompose the space of central charges into chambers. These are called $stablity \ chambers$ in \citep{L}. Now the statement of $Gaiotto-Moore-Neitzke \ 2d-4d \ wall \ crossing \ formula$\footnote{The name $2d-4d$ is from physics and we will not explain its meaning in this paper.} is:

{\bf When one changes the central charge $Z$, $A(\forall)$ is always constant as long as no BPS ($\mathcal{S}$ or $\mathcal{K}$) rays cross the boundary of $\forall$.}

This wall crossing formula becomes highly nontrivial algebraic identities in concrete examples, see \citep{GMN1}\citep{GMN2} for some examples. In particular it means that BPS numbers may change. The wall crossing formula is very powerful. In fact once the product is given and a stability chamber is chosen (hence an order of rays is given) the factorization into factors is determined (\citep{GMN2}). So if we know all nonzero BPS numbers of one stability chamber we can use the wall crossing formula to determine BPS numbers for other stability chambers.

A special case of this wall crossing formula is very important. Let us assume that $\mathcal{V},\Gamma_{i},\Gamma_{ij}$ are empty. Then we only have $\Omega(\gamma)$ as BPS numbers. In this case the wall crossing formula is called the $4d \ wall \ crossing \\ formula$ (see \citep{GMN2} for more details\footnote{Note that the conventions in \citep{GMN2} and \citep{GMN1} are slightly different due to the handling of $\sigma$ functions}). It is also known as the Kontsevich-Soibelman wall crossing formula and was first written down in \citep{KS}.

It is easy to see that it is consistent to consider only $X_{\gamma}$, $X_{\gamma_{ij}}$ and their transformations in the 2d-4d wall crossing formula. It only involves  $\Gamma, \Gamma_{ij}$ as lattices, $Z_{\gamma}, Z_{\gamma_{ij}}$ as central charges and $\Omega(\gamma), \omega(\gamma, \gamma_{ij}), \mu(\gamma_{ij})$ as BPS numbers. This kind of 2d-4d wall crossing formulas are called $restricted \ 2d-4d\  wall\  crossing \ formulas$ in this paper.

\subsection{ Hitchin's Moduli Spaces}

We will consider solutions of Hitchin's equations
\begin{equation}
\begin{array}{rr}F_{A}-\phi\wedge\phi &=0\\
d_{A}\phi=d_{A}\star\phi &=0
\end{array}
\end{equation} on a fixed Riemann surface $C$ of genus $g$. Here $A$ is a connection on a   bundle $E$ of rank 2 and degree 0 over $C$ whose gauge group is $G = SU(2)$, $\star$ is the Hodge star, $d_{A}:=d+A$, $\phi$ is an $ad(E)$-valued one form and $\wedge$ really means that we take the wedge of the one form part and the Lie bracket [ , ] of the bundle valued part. We will also use complex notations. In other words, we study a pair ($E$,$\varphi$) called a Higgs bundle or a Higgs pair. Here $E$ is a $holomorphic$-$G$ bundle and $\varphi$ is a $holomorphic$ one form valued in $ad(E)$. A Higgs bundle is obtained from a solution $(A, \phi)$ of equation (1) in the following way. The (0,1) part of $d_{A}$ denoted as $\bar{\partial}_{A}$ defines the holomorphic structure on $E$ and $\varphi$ is the (1,0) part of $i\phi$, $i\phi =\varphi+\bar{\varphi}$. $\varphi$ is also known as the Higgs field. In terms of Higgs bundles, the equivalent form of the equations are\begin{equation}\begin{array}{rr}F_{A}+[\varphi,\bar{\varphi}]&=0\\
\bar{\partial}_{A}\varphi &=0\\
\end{array}
\end{equation}
Following Seiberg-Witten and Gaiotto-Moore-Neitzke, we introduce an additional parameter $R$ and modify the  equations into \begin{equation}
\begin{array}{rr}F_{A}+R^{2}[\varphi,\bar{\varphi}]&=0\\
\bar{\partial}_{A}\varphi &=0\\
\end{array}
\end{equation}For the meaning and the significance of this parameter, see section 6.

 It is well known after Hitchin \citep{H1}\citep{H2} that the moduli space $\mathcal{M}$ of solutions of Hitchin's equations modulo the gauge equivalence is a  noncompact hyperkahler space. It is obtained by an infinite dimensional hyperkahler quotient construction of the moduli space. In fact the tangent space of the pair $(A,\phi)$ is an infinite dimensional affine space endowed with a natural flat hyperkahler metric (here the notations are the same as those of \citep{KW}) \begin{equation}ds^{2}=-{1\over 4\pi}\int_{C}|d^{2}z|\mathrm{Tr}(\delta A_{z}\otimes\delta A_{\bar{z}}+\delta A_{\bar{z}}\otimes\delta A_{z}+\delta \phi_{z}\otimes\delta \phi_{\bar{z}}+\delta \phi_{\bar{z}}\otimes\delta \phi_{z})\end{equation}where $A=A_{z}dz+A_{\bar{z}}d\bar{z}$ and $\phi=\phi_{z}dz+\phi_{\bar{z}}d\bar{z}$ for the holomorphic coordinate $z$ over $C$. $\delta$ denotes the tangent vectors. The group of gauge transformations acts on this flat hyperkahler space and the set of solutions of Hitchin's equations turns out to be the zero level set of associated three moment maps. Therefore $\mathcal{M}$ is obtained as a hyperkahler quotient.

 A hyperkahler space has a set of compatible complex structures parameterized by $\xi\in CP^{1}$ and generated by three independent complex structures $J_{1},J_{2},J_{3}$. The three independent complex structures satisfy the quaternion relations \begin{equation}J^{2}_{1}=J^{2}_{2}=J^{2}_{3}=J_{1}J_{2}J_{3}=-1\end{equation}The set of all compatible complex structures are given by \begin{equation}J_{\xi}:={i(-\xi+\bar{\xi})J_{1}-(\xi+\bar{\xi})J_{2}+(1-|\xi|^{2})J_{3}\over 1+|\xi|^{2}}\end{equation} where $\xi\in CP^{1}$  is called the $twistor \ parameter$. Let $\omega_{i}$ be the Kahler form in $J_{i}$. Then \begin{equation}\omega_{\xi}:={i(-\xi+\bar{\xi})\omega_{1}-(\xi+\bar{\xi})\omega_{2}+(1-|\xi|^{2})\omega_{3}\over 1+|\xi|^{2}}\end{equation}

  Define \begin{equation}\Omega(\xi)=-{i\over2\xi}\omega_{+}+\omega_{3}-{i\over2}\xi\omega_{-}\end{equation}
  \begin{equation}\omega_{\pm}:=\omega_{1}\pm i\omega_{2}\end{equation}Then $\Omega(\xi)$ is the holomorphic symplectic form in $J_{\xi}$. In particular$$\Omega_{1}=\omega_{3}-i\omega_{2}, \xi=i$$ $$\Omega_{2}=\omega_{3}+i\omega_{1}, \xi=-1$$ \begin{equation}\Omega_{3}=\omega_{1}+i\omega_{2}, \xi=0\end{equation}
 \begin{itemize}\item There are two opposite special complex structures which  are identified with $\pm J_{3}$ (i.e. $\xi$ is 0 or $\infty$). The moduli space $\mathcal{M}$ in $J_{3}$ is identified as the moduli space of semistable Higgs bundles.\item $\mathcal{M}$ is identified as the moduli space of  $SL(2,\mathbf{C})$ flat connections (the gauge group is $SL(2,\mathbf{C})$ because it is the complexification of $SU(2)$) when $\xi\neq 0,\infty$.  In fact, Hitchin's equations tell us that
the new connection $\mathcal{A}:={R\over\xi}\varphi+A+R\xi\bar{\varphi}$ is flat. \end{itemize}

There is a natural fibration from the moduli space to the space of quadratic differentials (we denote it by $B$) by taking the determinant of $\varphi$ (this is called Hitchin's map). It is called the $Hitchin's \ fibration$  and has the following properties\begin{itemize} \item The Hitchin's map denoted as $det$ is holomorphic, surjective and proper with respect to the complex structures in which the moduli space $\mathcal{M}$ is the moduli space of semistable Higgs pairs. \item Fibers of this map have nice geometric meanings. Define a curve in the total space of the canonical bundle of $C$ by the characteristic polynomial of $\varphi$
\begin{equation}det(x-\varphi)=0\end{equation}Note that the trace of $\varphi$ is zero. So the equation is $x^{2}+det\ \varphi=0$.  This curve is called the $spectral \ curve$ or the Seiberg-Witten curve $S$ and there is one such curve associated to each element of $B$. As an abelian variety, the fiber above $u\in B$ is  the Prym variety $\tilde{J}(S_{u})$ of the projection $S_{u}\rightarrow C$ (denote  the Jacobian  of $S_{u}$ by $J(S_{u})$. The Prym variety is the kernel of the natural map  $J(S_{u})\rightarrow J(C)$ induced by $S_{u}\rightarrow C$). Therefore we have  realized the moduli space as a family of complex abelian varieties.\item The complex dimension of $\mathcal{M}$ is $6g-6$. In particular we have to assume that $g> 1$. However later we will allow singularities of solutions of Hitchin's equations and this restriction will be removed then. \end{itemize}

\begin{dfn}\citep{GMN2} Let $S_{u}$ be a spectral curve of a Hitchin's moduli space $\mathcal{M}$. Let $\bar{S}_{u}$ be the compact Riemann surface obtained by filling in the punctures of $S_{u}$. We consider the odd part of $H_{1}(S_{u},\mathbf{Z})$. Here odd means that the cycle is invariant under the combined operation of exchanging the two sheets  and reversing the orientation. They fit into a local system over the nonsingular part of $B$. The local system is called the charge lattice and  degenerates at the singular locus of $B$ where some cycles become vanishing cycles. A charge is a section of the charge lattice local system. So locally by choosing an trivialization of the local system a charge is just an element of the associated integral lattice and  it has monodromies. The charge lattice is denoted as $\Gamma$ and is endowed with the antisymmetric intersection pairing of integral one cycles. The gauge charge lattice denoted by $\Gamma_{gau}$ is defined to be the local system of odd parts of $H_{1}(\bar{S_{u}},\mathbf{Z})$ together with the intersection paring. The flavor charge lattice $\Gamma_{flavor}$ is the radical of the intersection paring in $\Gamma$. It consists of integral combinations of loops around the punctures. Note these lattices  fit into an exact sequence $$0\rightarrow\Gamma_{flavor}\rightarrow\Gamma\rightarrow\Gamma_{gau}\rightarrow 0$$\end{dfn}

The genus of $C$ is $g$ while the genus of  $\bar{S_{u}}$ is $4g-3$. So the rank of $\Gamma_{gau}$ is $6g-6$ which matches the dimension of $B$.  Let $(A_{i},B_{i}), 1\leq i\leq3g-3$ be a symplectic basis of the gauge charge lattice. We can define  pairs of (half of them are redundant) holomorphic coordinates  by period maps\begin{equation}a_{i}(u):={1\over\pi}\int_{A_{i}}\lambda, \ \ a^{D}_{i}(u):={1\over\pi}\int_{B_{i}}\lambda\end{equation}We also define \begin{equation}\tau_{ij}(u):={da^{D}_{j}(u)/du\over da_{i}(u)/du}\end{equation} Note that $a_{i}(u)$ and $a^{D}_{i}(u)$ depend on
the holomorphic coordinate  on $B$    denoted by $u$ which is  defined by the value of the Hitchin's map $$u=det\ \varphi$$

\begin{dfn}The central charge $Z_{\gamma}$ for a charge  $\gamma$ is defined by \begin{equation}Z_{\gamma}={1\over\pi}\int_{\gamma}\lambda\end{equation}The central charge depends on $u$ since $\lambda$ does.\end{dfn}

The above theory holds for  more general types of Hitchin moduli spaces. We allow possibly singular solutions of Hitchin's equations. See \citep{GW} for a fraction of the huge literature. In this paper we only deal with regular singularities ($\varphi$ has only order one poles).

\subsection{4d Wall Crossing Formula and  Hyperkahler Metrics  on Hitchin's Moduli Spaces}

Let $\mathcal{M}$ be our Hitchin's moduli space. We view it as a moduli space of $SL(2,\mathbf{C})$ flat connections. Let us assume that there are only regular singularities $P_{i}, 1\leq i \leq l$ with regular semisimple residues. We assume $l\geq 1 $ in general and $l> 3$ if $g=0$.
\begin{dfn}\citep{GMN2} Choose a triangulation  of the Riemann surface $C$ with all vertices at singularities. Let $M_{i}$ be the clockwise monodromy of flat sections around $P_{i}$. Define a decoration at $P_{i}$ to be a choice of one of the two flat eigenlines of $M_{i}$. Denote such a decorated triangulation by $T$. For an edge $E$ of $T$, we consider the two triangles bounding $E$ making up a quadrilateral $Q_{E}$ with four vertices $P_{i},1 \leq i \leq 4$ in the counterclockwise order and $E$ connecting 1 and 3. Define the Fock-Goncharov  coordinate $\mathcal{X}_{E}^{T}$ by$$\mathcal{X}_{E}^{T}:=-{(s_{1}\wedge s_{2})(s_{3}\wedge s_{4})\over(s_{2}\wedge s_{3})(s_{4}\wedge s_{1})}$$where $s_{i}$ is an element of the one dimensional  decoration at $P_{i}$ (so it is defined up to a scaling and our definition is invariant under this scaling). Since $Q_{E}$ is simply connected, $s_{i}$ can be chosen to be single-valued in it and the four eigensections are evaluated at a common point $P_{*}$ inside the quadrilateral. The value is independent of the choice of the evaluation point because it is the $SL(2, \mathrm{C})$ invariant cross ratio.\end{dfn}

$\mathcal{X}_{E}^{T}$ is well defined on the Zariski open set which is the complement of the locus defined by the vanishing of the denominator. It is not hard to show that it is a holomorphic coordinate on this open subset. The set of all such functions where $E$ runs over all edges of a fixed decorated triangulation is a complete set of coordinates. Moreover outside the codimension one locus where either the numerator or the denominator is zero $\mathcal{X}_{E}^{T}$ is nonzero. So we have a set of locally defined $\mathbf{C}^{\times}$ valued functions. The dimension of the moduli space is \begin{equation}dim \mathcal{M}=6g-6+2l\end{equation}The number of edges are (see \citep{GMN2}) \begin{equation}\#E=6g-6+3l\end{equation} For each vertex (regular singularity) there is a relation between Fock-Goncharov coordinates given by taking the product of all that are associated to edges meeting at that vertex. So the number of independent ones is the same as $dim\mathcal{M}$.

We allow degenerate triangulations. Here a degenerate triangulation means that  two edges in a triangle are identified. So we have a double vertex\footnote{Actually we also allow the even more degenerate situation where we have triple vertices. The Fock-Goncharov formalism still works.} and the edge connecting the double vertex to itself is a loop while the double edge is an edge connecting a point on the loop (the double vertex) to another vertex $P$. To define Fock-Goncharov coordinates in this situation, we take a cover ramified at $P$ such that after taking the pull-back the triangulation is non-degenerate. Such a cover always exists and is non-unique but our definition does not depend on the choice of the cover. We pull back everything and define the Fock-Goncharov  coordinate for the degenerate edge $E$ to be the ordinary Fock-Goncharov  coordinate $\mathcal{X}_{\tilde{E}}^{T}$ on the cover where $\tilde{E}$ is any choice of the pre-images of $E$ and the definition does not depend on this choice.  The dimensional count of Fock-Goncharov coordinates is still valid for degenerate decorated triangulations.

Let us consider the Hitchin's fibration.   The determinant of $\varphi$ is a quadratic differential $-\lambda^{2}$ well defined on $C$ and therefore we take the holomorphic coordinate $u$ of the base $B$ of the Hitchin' fibration (recall that $B$ can be identified as  the space of quadratic differentials) to be $u= -\lambda^{2}$. $\lambda^{2}$ has order two poles and generically has only simple zeroes. $\lambda$ is a one form  on the Riemann surface $C$ defined up to a sign but is a single valued one form on the spectral curve $S$ which is a double cover of $C$.

Fix an angular parameter $\vartheta$ and consider the foliation given by trajectories of $\lambda^{2}$ with phase $\vartheta$.
\begin{dfn}A trajectory of $\lambda^{2}$ with phase $\vartheta$ is a curve whose tangent vector $\partial_{t}$ satisfies $$\langle\lambda, \partial_{t}\rangle\in e^{i\vartheta}\mathrm{R}^{\times}$$ everywhere on the curve.\end{dfn}
\begin{note}A trajectory is called a WKB curve in \citep{GMN2}.\end{note}
There is an extensive theory of foliations given by meromorphic quadratic differentials  and the local behaviors near singularities and zeroes as well as global behaviors are known. The standard reference is Strebel' book \citep{St}. Let us summarize the results we need.

Near a point which is neither a zero nor a pole of the quadratic differential, we can straighten the foliation by choosing local coordinate $w:=\int\lambda$. Locally near an order $n$ zero, we can choose a local parameter $\zeta$ such that $\lambda^{2}$ has the representation \begin{equation}\lambda^{2}=({n+2\over 2})^{2}\zeta^{n}d\zeta^{2}\end{equation}The full angle $0\leq\arg\zeta\leq2\pi$ is divided into $n+2$ equal sectors. In particular for a simple zero (i.e. order one zero) the foliation develops three asymptotic directions surrounding and going away from the zero.

Since we only have regular singularities for the Hitchin's equations the order of poles of $\lambda^{2}$ is two. It is shown that in this case it has a local representation of the form $$\lambda^{2}={a\over\zeta^{2}}d\zeta^{2}$$ and the trajectories near the pole is either logarithmic spirals approaching the pole or radii approaching the pole or closed circles around the pole.

Globally a trajectory belongs to one of the following cases.\begin{itemize}\item A generic trajectory. It is asymptotic in both directions to singular points. Generic trajectories arise in one dimensional families.\item A separating trajectory.  It is asymptotic in one direction to a simple zero and in the other direction to a singular point. Separating trajectories separate families of generic trajectories. \item A finite trajectory. It is asymptotic in both directions to a simple zero (both directions could go to the same zero) or  is closed.\item A divergent trajectory. It is neither closed nor approaches to a limit in one or both directions.\end{itemize}

For a generic $\vartheta$,  finite trajectories are absent and in that case Gaiotto, Moore and Neitzke showed the absence of divergent ones in our setting. We will assume the absence of finite trajectories for now to get decorated triangulations and later we will see they are the source of BPS numbers.

Following Gaiotto, Moore and Neitzke, we define a decorated triangulation called WKB triangulation in the following way. \begin{dfn}We take a generic $u=-\lambda^{2}$ such that it has only simple zeroes (note that this is the generic case, the nongeneric ones are codimensional two in $B$). We choose one element from every family of generic trajectories separated by separating trajectories. They make  an ordinary triangulation of the Riemann surface. The choices of the representatives are unimportant because a triangulation is only meant to be defined up to isotopy.

There might be  a generic trajectory approaching the same singularity along both directions in which case we get a degenerate triangulation.

Near each singularity, there are two independent eigen flat sections. It is shown that one of them is exponentially small along trajectories going to the singularity while the other is exponentially large. We pick the small flat section as the decoration at the singularity. These decorations together with the triangulation define a decorated triangulation $T_{WKB}(\vartheta, \lambda^{2})$ called a WKB triangulation and therefore a set of Fock-Goncharov coordinates $\mathcal{X}_{E}^{T_{WKB}(\vartheta, \lambda^{2})}$.\end{dfn}

Fock-Goncharov coordinates are labeled by edges. We want them to be labeled by charges. It is easy to see that every triangle in a WKB triangulation contains exactly one simple zero. Let $E$ be the edge labeling the Fock-Goncharov coordinate $\mathcal{X}_{E}^{T}$. We choose an oriented simple loop inside $Q_{E}$ surrounding the two zeroes in the two adjacent triangles and define the associated charge $\gamma_{E}$ to be the lift of the loop to the spectral curve $S$ which is a double cover of the underlying Riemann surface. Ambiguity of the sign of the cycle induced by ambiguities  of choosing orientations and one of the two sheets can be canonically fixed in the following way. Note  that $\lambda$ is a single-valued one form over $S$. We require that the positively oriented tangent vector $\partial_{t}$ of the lift of $E$ to $S$ denoted as $\hat{E}$ satisfies \begin{equation}e^{-i\vartheta}\langle\lambda,\partial_{t}\rangle>0\end{equation} The sign of the cycle $\gamma_{E}$ is fixed by $\langle\gamma_{E},\hat{E}\rangle=1$. So we can replace the labeling by $E$ by labeling by $\gamma_{E}$ (also denoted simply as $\gamma$ later in this paper). This operation respects the integral antisymmetric pairings, i.e. $\langle\gamma_{E},\gamma_{E^{'}}\rangle=\langle E,E^{'}\rangle$. It is easy to generalize to degenerate edges.

It is shown in \citep{GMN2} that one can generate the charge lattice $\Gamma$  by  cycles associated to edges. We write \begin{equation}\mathcal{Y}_{\gamma_{E}}=-\mathcal{X}_{\gamma_{E}}\end{equation}We then extend the definition of Fock-Goncharov coordinates to the whole lattice by the multiplicative relation \begin{equation}\mathcal{Y}_{\gamma_{E}}\mathcal{Y}_{\gamma_{E^{'}}}=(-1)^{\langle\gamma_{E},\gamma_{E^{'}}\rangle}\mathcal{X}_{\gamma_{E}+\gamma_{E^{'}}}\end{equation}

Here $\mathcal{Y}_{\gamma}$ is actually the short notation for $\mathcal{Y}^{\vartheta,u}_{\gamma}(x;\xi)$ where we have restored all the variables involved in the definition. $x\in \mathcal{M}$ and $\xi \in H_{\vartheta}$.\begin{equation}H_{\vartheta}:=\{\xi\mid \vartheta-\pi/2<arg \xi<\vartheta+\pi/2\}\end{equation}The dependence on $x,\vartheta,u$ is clear. It also depends on $\xi$ because the flat connection depends on it and here we restrict $\xi$ to $H_{\vartheta}$ for a given $\vartheta$.

We can vary the phase $\vartheta$. Although for a generic $\vartheta$ there are no finite trajectories they do appear for exceptional values of phases of $\vartheta$. If we label $\vartheta$ by rays on the complex plane then when $\vartheta$ crosses countably many exceptional rays, the decorated triangulation defined up to isotopy will change and we encounter discontinuous  transformations of Fock-Goncharov coordinates.

These nontrivial transformations of WKB triangulations are classified into three types: \begin{itemize}\item A flip. An edge is flipped when an exceptional ray of $\vartheta$ is crossed. By flip at an edge we mean that we replace $E=E_{13}$ which connects vertices 1 and 3 in the quadrilateral containing $E$ as a diagonal edge by $E^{'}=E_{24}$ connecting vertices 2 and 4 and obtain a new triangulation $T^{'}$. As $\vartheta$ goes to the exceptional ray the flipped edge degenerated to a finite trajectory connecting two simple zeroes. There is only one such finite trajectory for that exceptional value of $\vartheta$.\item A pop for a degenerate triangle. We will not describe it as it leads to trivial transformations of Fock-Goncharov coordinates. \item Infinitely many flips leading to a limit configuration. In this case, the phase first passes infinitely many rays (corresponding to  flips) to reach a special configuration with  closed trajectories  around the pole inside an annular region bounded by two finite boundary trajectories. Each boundary trajectory's  both directions go to a same simple zero (but the simple zeroes are different for the two trajectories). We can take the limit of Fock-Goncharov coordinates. On the other hand, if we start from the other side of the special ray and approach it from the other direction we would also pass infinitely many rays and get another set of limit Fock-Goncharov coordinates. The two limits are differed by a nontrivial transformation   called a juggle.  \end{itemize}

Moreover, we can vary $u$. The exceptional phases clearly depend on $u$. The order of exceptional phases does not change if we take small enough perturbations of $u$. But it will change if we let $u$ cross codimension one locus of $B$ which corresponds to the coincidence of two exceptional phases.

\begin{thm}\citep{GMN2} Suppose the real parameter $R$ in the definition of Hitchin's equations is large enough. The exceptional phases $\vartheta$ mentioned above are phases of BPS $\mathcal{K}$-rays and vice versa. One can assign BPS numbers $\Omega(\gamma)$ to them ($\gamma$ is the corresponding charge of the BPS $\mathcal{K}$-ray) and define appropriate twisting functions such that
\begin{enumerate}\item For each $\xi \in\mathbf{C}^{\times}$, $\mathcal{Y}_{\gamma}$ is valued in $C^{\times}$ and is holomorphic in the complex structure $J_{\xi}$.\item \begin{equation}\mathcal{Y}^{\vartheta,u}_{\gamma}(x;\xi)=\bar{\mathcal{Y}}^{\vartheta+\pi,u}_{-\gamma}(x;-1/\bar{\xi})\end{equation}where the right hand side is the complex conjugate of $\mathcal{Y}^{\vartheta+\pi,u}_{-\gamma}(x;-1/\bar{\xi})$.\item $\mathcal{Y}_{\gamma}$ is holomorphic in $\xi$ for $\xi\in H_{\vartheta}$.
\item \begin{equation}\lim_{\xi\rightarrow 0}\mathcal{Y}_{\gamma}\exp[-\xi^{-1}\pi RZ_{\gamma}]\end{equation} exists.
\item The discontinuous transformations of $\{\mathcal{Y}_{\gamma}\}$ are given by $\mathcal{K}$-factors.\item The 4d wall crossing formulas holds\footnote{We identify $\mathcal{Y}_{a}$ as the formal variable $X_{a}$.} (this statement makes sense because according to the definition of $Z$, $B$ can be identified with the space of central charges $Z$)\footnote{Since the explicit form of the 4d wall crossing formula involves $\mathcal{K}$-factors from now on when we say some data satisfies the 4d (2d-4d) wall crossing formula we always mean that the discountinuities are $\mathcal{K}$-factors (or $\mathcal{K}$-factors and $\mathcal{S}$-factors). and the 4d (2d-4d) wall crossing formula  holds.}\end{enumerate}
\end{thm}
The details can be found in \citep{GMN2}. Roughly speaking when finite trajectories appear  one can associate a charge $\gamma$ and a nonzero $\Omega(\gamma)$ to it such that the transformations of $\{\mathcal{Y}_{\gamma}\}$ are given by $\mathcal{K}$-factors. For example for a flip $\gamma$ is the charge corresponding to the flipped edge in the above charge labeling construction. The BPS numbers are \begin{equation}\Omega(\gamma)=1,\ \Omega(n\gamma)=0, n\neq 1\end{equation}Note that the 4d wall crossing formula  in \citep{GMN2} (called the wall crossing formula there) looks slightly different from the one in \citep{GMN1}. To recover the one in \citep{GMN1} we need to transform the so-called quadratic refinements of \citep{GMN2} to the twisting functions in \citep{GMN1} by the recipe of section 7 of \citep{GMN1}.

Let $\mathcal{M}$ be an $SU(2)$ Hitchin's moduli space with only regular singularities $P_{i}, 1\leq i \leq l$ with regular semisimple residues. We assume $l\geq 1 $ in general and $l> 3$ if $g=0$. Recall that the definition of $\mathcal{M}$ depends on $R$.
\begin{thm}For a large enough $R$, Fock-Goncharov coordinates labeled by charges define a hyperkahler structure on the Hitchin's moduli space $\mathcal{M}$ via the twistor construction.\end{thm}

We are not going to explain what the twistor construction is. But we want to point out that in this construction Fock-Goncharov coordinates are not directly viewed as coordinates on $\mathcal{M}$. We use them together with their dependence on $\xi$ to produce coordinates on $\mathcal{M}\times CP^{1}$ which is the starting point of the so-called Gaiotto-Moore-Neitzke ansatz for a twistor construction. From this perspective the 4d wall crossing formula is a consistency condition. Roughly speaking the special form of discontinuous jumps guarantees that the twistor construction is well defined even if Fock-Goncharov coordinates are discontinuous. However this well-definedness condition itself inevitably depends on a choice of all charges with nonzero BPS numbers which depends on a choice of stability chamber and hence has a potential inconsistency. The 4d wall crossing formula then closes this gap\footnote{A thorough discussion can be found in \citep{L}}. Later when we discuss an analogous twistor construction of hyperholomorphic connections we will meet the 2d-4d wall crossing formula which plays the same role there.

The proofs of the above two theorems are long and complicated. They are given by Gaiotto-Moore-Neitzke and are contained in \citep{GMN2} and \citep{GMN3}. An exposition for mathematicians can be found in \citep{L}(which also fills some minor gaps).

It is natural to conjecture that \begin{con}The hyperkahler metric constructed above is isometric to the hyperkahler quotient metric of Hitchin. \end{con}

\subsection{Riemann-Hilbert Problems}

There is another way to understand the above data which was introduced in \citep{GMN3}. From $\mathcal{Y}^{\vartheta}_{\gamma}(\xi)$ (we only restore the relevant variables) one divides the $\xi$-plane into slivers bounded by BPS $\mathcal{K}$-rays and then defines a new function $\mathcal{Y}^{RH}_{\gamma}(\xi)$ to agree with $\mathcal{Y}^{\vartheta}_{\gamma}(\xi)$ in the sliver containing the ray $e^{i\vartheta}\mathbf{R}^{+}$. $\mathcal{Y}^{RH}_{\gamma}(\xi)$ then are piecewise holomorphic and has discontinuous jumps given by the $\mathcal{K}$-factors. We can recover $\mathcal{Y}^{\vartheta}_{\gamma}(\xi)$ from $\mathcal{Y}^{RH}_{\gamma}(\xi)$ by taking analytic continuations  in $\xi$ from $e^{i\vartheta}\mathbf{R}^{+}$.

The function $\mathcal{Y}^{RH}_{\gamma}(\xi)$ can be considered as a solution of the following Riemann-Hilbert type problem:

Find a piecewise holomorphic function with discontinuities on the $\xi$-plane given by $\mathcal{K}$-rays and $\mathcal{K}$-factors and asymptotic behaviors as $\xi\rightarrow 0, \infty$ given by  the asymptotics of the Fock-Goncharov coordinates in theorem 1.1.

Gaiotto, Moore and Neitzke write down a formulation of the above Riemann-Hilbert problem in terms of integral equations. Suppose we can
solve the following integral equations (say by some iteration scheme) \begin{equation}\mathcal{Y}_{\gamma}(\xi)=\mathcal{Y}_{\gamma}(\xi)^{sf}\exp({1\over 4\pi i}\Omega(\gamma^{'})\langle\gamma^{'}, \gamma\rangle\sum_{l_{\gamma^{'}}}\int_{l_{\gamma^{'}}}{d\xi^{'}\over \xi^{'}}{\xi^{'}+\xi\over \xi^{'}-\xi}\log(1-\mathcal{Y}_{\gamma^{'}}(\xi^{'}))\end{equation}then we have solved the Riemann-Hilbert problem. Here the sum runs over all BPS $\mathcal{K}$-rays $l_{\gamma^{'}}$.
$\mathcal{Y}^{sf}_{\gamma}$ is defined below and it guarantees the expected asymptotic behaviors.

The following conjecture essentially says the solution of the Riemann-Hilbert problem is unique.
\begin{con}$\mathcal{Y}^{RH}_{\gamma}(\xi)$ satisfies equation (34) if we replace $\mathcal{Y}$ by $\mathcal{Y}^{RH}$.\end{con}

\subsection{B Fields}
We need some knowledge about universal bundles, B fields and twisted bundles. The following exposition follows \citep{H3} and \citep{KW}.

Recall that the $SU(2)$  Hitchin's moduli space $\mathcal{M}$ carries a hyperkahler quotient metric obtained from the hyperkahler quotient construction. Suppose $\mathcal{M}$ is smooth. we want to have a universal bundle over $\mathcal{M}\times C$ or at least over $\mathcal{M}\times (C\setminus singularities)$ in the sense that  the fiber over a point of $\mathcal{M}$ is the solution (up to gauge equivalence) of Hitchin's equations on $C$ corresponding to this point. It turns out that we can not have such a universal bundle due to nontrivial automorphisms but we can have a twisted universal bundle.

Let $z$ be a point of $C$ which is not a singularity of solutions of Hitchin's equations. Let $\mathcal{G}_{z}$ be the subgroup of the gauge group consisting of gauge transformations which restrict to the identity at $z$. To get $\mathcal{M}$ we can first take the hyperkahler quotient as in section 2.2 but with gauge group $\mathcal{G}_{z}$. This produces a space denoted by $\mathcal{N}_{z}$. If we take the hyperkahler quotient of it by $SU(2)$ then we get $\mathcal{M}$. Here the center of $SU(2)$ is $\mathbf{Z}_{2}$ which acts trivially on $\mathcal{N}_{z}$ while $PSU(2)=SO(3)=SU(2)/\mathbf{Z}_{2}$ acts freely. The hyperkahler quotient of $\mathcal{N}_{z}$ by  $PSU(2)$ gives us a principal $PSU(2)$ bundle (which is the zero locus of the hyperkahler moment map)  over $\mathcal{M}$.

We are going to explain that this bundle (which is denoted by $\mathcal{U}_{z}$) is a  twisted bundle over $\mathcal{M}$ twisted by a B field over $\mathcal{M}$. We shall use $V_{z}$ to denote the associated twisted vector bundle (see below) of $\mathcal{U}_{z}$ associated to the defining representation.

To define B fields we follow Hitchin \citep{H3}.

\begin{dfn}Let $X$ be a manifold and $\{U_{\alpha}\}$ be an open covering. A $U(1)$ gerbe (which will be just called a gerbe in this paper) is a data consisting of a map $$g_{\alpha\beta\gamma}:U_{\alpha}\cap U_{\beta}\cap U_{\gamma}\rightarrow U(1)$$ for each threefold intersection with $$g_{\alpha\beta\gamma}=g_{\beta\alpha\gamma}^{-1}=g_{\alpha\gamma\beta}^{-1}=g_{\gamma\beta\alpha}^{-1}$$and satisfying the cocycle condition \begin{equation}g_{\beta\gamma\delta}g_{\alpha\gamma\delta}^{-1}g_{\alpha\beta\delta}g_{\alpha\beta\gamma}^{-1}=1\ on\ U_{\alpha}\cap U_{\beta}\cap U_{\gamma}\cap U_{\delta}\end{equation}
The cocycle defines a class in $H^{2}(X,C^{\infty}(U(1)))\simeq H^{3}(X,\mathbf{Z})$. Similarly one can define a $\mathbf{Z}_{n}$ gerbe.

A specific representation of the cocycle $g_{\alpha\beta\gamma}$ as a coboundary is called a trivialization of the gerbe. The difference of two trivializations is a line bundle.

A connection of on a gerbe is the data consisting of 1-forms $\{A_{\alpha\beta}\}$ on twofold intersections, 2-forms $\{F_{\alpha}\}$ on $\{U_{\alpha}\}$ and a closed global 3-form $G$ satisfying $$G|_{U_{\alpha}}=dF_{\alpha}$$ $$F_{\beta}-F_{\alpha}=dA_{\alpha\beta}$$
\begin{equation}iA_{\alpha\beta}+iA_{\beta\gamma}+iA_{\gamma\alpha}=g_{\alpha\beta\gamma}^{-1}dg_{\alpha\beta\gamma}\end{equation}
$G$ is called the curvature. A connection on a gerbe is flat if its curvature is zero. A flat gerbe is called a  B field.\end{dfn}

For a  B field one can define its holonomy (see \citep{H3}) which is a  class in $H^{2}(X, U(1))$ just like the holonomy of a flat connection on a unitary bundle defines a class in $H^{1}(X, U(1))$.  Now suppose the holonomy is trivial in which case we say that the B field is isomorphic to a trivial B field. One can define a flat trivialization (see \citep{H3}) of a B field isomorphic to a trivial B field. It is the analog of a covariantly constant trivialization of a flat line bundle. The difference of two flat trivializations is a flat line bundle. We can also define the tensor product of B-fields. A more general formalism for gerbes and B-fields uses the notion of sheaf of categories. We will try to avoid this language in this paper.\\

Back to our universal bundle problem. Locally in $\mathcal{M}$ in the product $\mathcal{M}\times C$ there is no obstruction to the existence of the universal $SU(2)$ bundle (i.e. the restriction to a point of $\mathcal{M}$ gives us the corresponding $SU(2)$ Higgs bundle over $C$). The obstruction to the existence of the global universal bundle is given by the cocycle defined by \begin{equation}f_{\alpha\beta\gamma}:=\Theta_{\alpha\beta}\Theta_{\beta\gamma}\Theta_{\gamma\alpha}\end{equation}where $\{\Theta_{\alpha\beta}\}$ are transition functions. Note that $f_{\alpha\beta\gamma}$ is valued in the automorphism group which is the center $\mathbf{Z}_{2}$. This defines a class in $H^{2}(\mathcal{M},\mathbf{Z}_{2})$ corresponding to a $\mathbf{Z}_{2}$ gerbe. Let $\rho$ be the defining representation of $SU(2)$ then the center acts by scalar multiplications. $\rho(f_{\alpha\beta\gamma})$ then is a cocycle defining a class in $H^{2}(\mathcal{M}, C^{\infty}(U(1)))$. The data consisting of such transition functions  is called a twisted vector bundle, twisted by $\rho(f_{\alpha\beta\gamma})$. So after specifying a class in $H^{2}(\mathcal{M}, C^{\infty}(U(1)))$ we have a twisted vector bundle.

The bundle  $\mathcal{U}_{z}$ is clearly a twisted principal bundle for the gauge group $SU(2)$ while $V_{z}$ is a twisted vector bundle for the gauge group $SU(2)$. Let $\zeta$ and $\rho(\zeta)$ denote the associated classes in $H^{2}(\mathcal{M},\mathbf{Z}_{2})$ and $H^{2}(\mathcal{M}, C^{\infty}(U(1)))$ respectively. The same analysis works for the Higgs bundle picture. In fact the obstruction of lifting from $PSU(2)$ to $SU(2)$ for Higgs bundles is given by the bundle part of the Higgs pair and hence produces the same twisting class.

The above discussion can be extended to the case where  Hitchin's moduli spaces $\mathcal{M}$ have orbifold singularities obtained by taking  global quotients by the actions of  finite groups, see page 5 of \citep{HT}. The moduli spaces we study in this paper have only such singularities.

\subsection{2d-4d Wall Crossing Formula and Hyperholomorphic Connections on Hitchin's Moduli Spaces}

Gaiotto, Moore and Neitzke have developed a new construction of hyperholomorphic connections on twisted bundles $V_{z}$.
\begin{dfn}A connection over a (possibly twisted) vector bundle over a hyperkahler space is hyperholomorphic if its curvature is type $(1,1)$ in all compatible complex structures.\end{dfn}
By treating the $(0,1)$ part of a hyperholomorphic connection as the $\bar{\partial}$ operator we can give the bundle a holomorphic structure for each compatible complex structure. There is a tautological hyperholomorphic connection on $V_{z}$ according to \citep{N}, but  Gaiotto, Moore and Neitzke propose a construction of much more general ones from the quadratic differential foliations used in the construction of Fock-Goncharov coordinates.

Let $\mathcal{M}$ be the $SU(2)$ Hitchin's moduli space with regular singularities as in section 2.3. We use the same definitions of section 2.3. So we have defined  $\Gamma$, $Z:\Gamma\rightarrow \mathbf{C}$, $\mathcal{Y}_{\gamma}$ and $\Omega(\gamma)$. Let $z \in C$ be a point which is not one of the regular singularities. Define the finite set $\mathcal{V}_{z}$ to be the set of the two preimages in the double cover $\Sigma\rightarrow C$ where $\Sigma$ is the spectral curve. We use $\{x_{i}\},\ i=1,2$ to denote its elements. The $\Gamma$-torsor $\Gamma_{ij}$ is defined to be the $\Gamma$-torsor of relative homology classes of open paths from $x_{i}$ to $x_{j}$. So we just pick   a path from $x_{i}$ to $x_{j}$ in $\Sigma$ and then collect all translates by $\Gamma$. The central charges $Z_{\gamma_{ij}}$ are defined by
\begin{equation}Z_{\gamma_{ij}}:={1\over\pi}\int_{\gamma_{ij}}\lambda\end{equation}
We also define $\omega(\gamma,\gamma_{ij})$. When $\gamma$ is a charge associated to a flip, let $\hat{\gamma}$ be the lift of the finite trajectory to $\Sigma$ described in section 2.3 (so $\hat{\gamma}$ is a representative of the class $\gamma$) and define \begin{equation}\omega(\gamma,\gamma_{ij}):=\Omega(\gamma)\langle\hat{\gamma}, \gamma_{ij}\rangle\end{equation}When $\gamma$ is a charge associated to a juggle, the definition is more subtle and can be found in section 7 of \citep{GMN2}. The definition of the twisting functions can also be found there.

$\mu(\gamma_{ij})$ is defined to be the number of so-called finite open WKB networks. That means  finite trajectories or separating trajectories which contain $z$. The part of the lift  of a finite open WKB network between preimages of $z$ gives us an element $\gamma_{ij}$ of $\Gamma_{ij}$. For our $\mathcal{M}$ there is at most one finite open WKB network for each $\gamma_{ij}$. The number $\mu_{\gamma_{ij}}$ is defined to be either $1$ or $0$ depending on whether the finite open WKB network exists or not.

Just like we define $\mathcal{Y}_{\gamma}$ which are identified as $X_{\gamma}$ in the 4d wall crossing formulas we want to define $\mathcal{Y}_{\gamma_{ij}}$ which are identified as $X_{\gamma_{ij}}$ in the 2d-4d wall crossing formula. Consider a WKB triangulation and let $S$ be the sector containing $z$. Here a sector is a triangle bounded by an edge of the WKB triangulation and two separating WKB curves. So two vertices are regular singularities while the third one is a simple zero of the quadratic differential. Define \begin{equation}s_{i,S}:=(s_{b}\wedge s_{c})s_{a}\end{equation} where $a$ is the vertex of the sector $S$ such that the lift of the WKB trajectory passing though $x_{i}$ goes to it. So it must be one of the two singularities. The vertices $abc$ go around the sector counterclockwise. We define $s_{j,S}$ similarly.

 Consider the twisted bundle $V_{z}$ associated to $z$. Note that $s_{i,S}(z)$ and $s_{j,S}(z)$ are sections of this bundle. We define an endomorphism $\mathcal{Y}_{\gamma_{ij},S}$ of it by \begin{equation}\mathcal{Y}_{\gamma_{ij},S}(s_{i,S})=0,\ \ \mathcal{Y}_{\gamma_{ij},S}(s_{j,S})=\nu_{i,S}s_{i,S}\end{equation}where $\nu_{i,S}=1$ if the lifted WKB trajectory through $x_{i}$ goes around the sector counterclockwise and $\nu_{i,S}=-1$ if otherwise. We also define \begin{equation}\mathcal{Y}_{\gamma_{ij}+\gamma}=\sigma(\gamma)\mathcal{Y}_{\gamma}\mathcal{Y}_{\gamma_{ij}}\end{equation}
where $\sigma(\gamma)$ is the canonical quadratic refinement defined in section 7 of \citep{GMN2}. For practical purpose we only need to know that it is $-1$ if $\gamma$ is associated to a flip and is $1$ if $\gamma$ is associated to a juggle.

The data consisting of $\mathcal{V}_{z}$, $\Gamma$, $\{Z_{\gamma}\}$, $\{\mathcal{Y}_{\gamma}\}$, $\{\Omega(\gamma)\}$, $\{\Gamma_{ij}\}$, $\{Z_{\gamma_{ij}}\}$, $\{\omega(\gamma,\gamma_{ij})\}$, $\mathcal{Y}_{\gamma_{ij}}$ and twisting functions is enough for setting up the restricted wall crossing formula. Now we can state one of the main results in \citep{GMN1}:
\begin{thm}\citep{GMN1}The above data satisfies the restricted 2d-4d wall crossing formula.\end{thm}Details can be found in section 7 and appendix F of \citep{GMN1}.

We should point out that for 2d-4d wall crossing the central charges $Z_{\gamma_{ij}}$ depend on not only the quadratic differential but also the location of $z$. And we have not only $\mathcal{K}$-factors but also $\mathcal{S}$-factors because the change of $\vartheta$ can cause not only nontrivial changes of the underlying WKB triangulation but also  changes of the sector that contains $z$.

We now have endomorphisms of the bundle $V_{z}$ i.e.  sections of $Hom(V_{z}, V_{z})$.\footnote{In \citep{GMN1} there is a more general construction producing  sections of $Hom(V_{z}, V_{z^{'}})$ where $z^{'}$ is a different point of $C$. We will not need this more general construction.} We want to construct sections of $V_{z}$. To do that we will first try to define $\Gamma_{i}$. $\Gamma_{i}$ is supposed to be associated to one of the two preimages of $z$ and we must have the splitting of charges \begin{equation}\Gamma_{ij}=\Gamma_{i}-\Gamma_{j}\end{equation} and central charges \begin{equation}Z_{ij}=Z_{i}-Z_{j}\end{equation}Here $i, j=1,2$. There is an obvious obstruction (even before we describe how to define $\Gamma_{i}$) to the existence of a local system of torsors $\Gamma_{i}$ on the complement of the discriminant locus of $B$. In fact monodromies around the discriminant locus of $B$ can change the sheets of spectral curves and carry an element of $\Gamma_{i}$ to $\Gamma_{j}$. We can take a ramified double cover $B_{S}$ of $B$ such that the sheets of $B_{S}$ correspond to the sheets of  spectral curves. Then we go on and combine $\Gamma_{1}$ and $\Gamma_{2}$ to form a local system of torsors on the complement of the preimage of the discriminant of locus in $B_{S}$ \citep{GMN1}. We still use $\Gamma_{i}$ to denote this local system although now $i$ has included both sheets. However it turns out that there could be another obstruction. It is called an anomaly in \citep{GMN1}.

Pick a representative $\gamma_{12}\in\Gamma_{12}$. Gaiotto, Moore and Neitzke define \begin{equation}\Gamma_{1}:={1\over 2}\gamma_{12}+\Gamma\end{equation} and then define $Z_{\gamma_{1}}$ by linearity. One defines $\Gamma_{2}$ in the same way. It seems then we have got the required splitting properties. However since $\Gamma$ is a local system of lattices there are monodromies around loops in $B_{S}$. That means that even though we have taken care of the change of sheets by working over $B_{S}$ there is still an  ambiguity of $\gamma_{ij}$ given by shifting  by an element of $\Gamma$. Then we can only guarantee that the ambiguity of an element of $\Gamma_{i}$ is an element of ${1\over 2}\Gamma$. This fractional shift is the obstruction.

 Let $\theta_{\gamma^{e}_{I}}, \theta_{\gamma^{m}_{I}}$ be angular coordinates of the fibers of the Hitchin's fibration denoted by $\pi$ (in \citep{GMN1} they are denoted by $\theta_{e}^{I},\ \theta_{m,I}$). Here $\{\gamma^{e}_{I}, \gamma^{m}_{I}\}$ is a symplectic basis of the gauge charge lattice. We can identify the integral local system underlying the torus fibration $\pi$ away from singular fibers  as the dual of $\Gamma_{gau}$ ($\Gamma_{gau}$ is self dual) and then locally in $B^{*}$ (the complement of the discriminant locus of $B$) we can identify $\theta_{\gamma^{e}_{I}}, \theta_{\gamma^{m}_{I}}$  with local sections of the local system of gauge charges tensored with $2\pi\mathbf{R}/\mathbf{Z}$. We have \begin{equation}a_{I}(u):={1\over\pi}\int_{\gamma^{e}_{I}}\lambda, \ \ a^{D}_{I}(u):={1\over\pi}\int_{\gamma^{m}_{I}}\lambda\end{equation} \begin{equation}\tau(u)_{IJ}:={da^{D}_{J}(u)/du\over da_{I}(u)/du}\end{equation} $Z_{\gamma_{i}}$ defined above depends on $u$ and hence on $a$ if we use the holomorphic coordinates $a_{I}$ instead of $u$ to parameterize $B$. Define\begin{equation}t^{i}_{I}:=\eta_{I}^{i}+\tau_{IJ}\alpha_{J}^{i}:= {\partial Z_{\gamma_{{i}}}\over \partial a_{I}}\end{equation} where $\eta_{I}^{i},\alpha_{I}^{i}$ are real functions of $a_{I}$. They are determined  by $t^{i}_{I}, \tau_{IJ}$. Now consider the following locally defined one form on $\mathcal{M}$:\begin{equation}A_{i}^{sf}:=i\sum_{I}(\eta^{i}_{I}d\theta_{\gamma_{I}^{e}}+\alpha^{i}_{I}d\theta_{\gamma_{I}^{m}})\end{equation}

Let $\{\mathcal{U}_{\alpha}\}$ be an open covering of $B^{*}$. $\{\pi^{-1}(\mathcal{U}_{\alpha})\}$ is an open covering of the complement of singular fibers in $\mathcal{M}$. In fact we need to go to the cover $B^{*}_{S}$ and replace $\mathcal{M}$ by the fiber product $\mathcal{M}\times_{B}B_{S}$. We keep the  same notations  for objects and their pullbacks to $B_{S}$ or the fiber product.  Locally as an element of $\Gamma$-torsor $\gamma_{i}$ can be shifted by some $\tilde{\gamma}^{i}\in\Gamma$. This leads to a shift \begin{equation}A_{i}^{sf}\rightarrow A_{i}^{sf}+id\theta_{\tilde{\gamma}^{i}}\end{equation} Define transition functions on overlaps by $\{e^{i\theta_{\tilde{\gamma}^{i},\ \alpha\beta}}\}$ where the subscript $\alpha\beta$ indicates that we are in the twofold intersection $\mathcal{U}_{\alpha\beta}$. We then interpret $A_{i}^{sf}$ globally as a connection of the line bundle defined by these transition functions. However the fractional charge ambiguity (anomaly) means that this is actually a twisted line bundle. In fact the ${1\over 2}\Gamma$ anomaly implies that in $\mathcal{U}_{\alpha\beta}$ we have an ambiguity by a square root of unity for $e^{i\theta_{\tilde{\gamma}^{i},\ \alpha\beta}}$. If we make a choice in $\mathcal{U}_{\alpha\beta}$ the cocycle condition in threefold intersections will be violated by a square root of unity which gives us a class  in $H^{2}(\mathcal{M},C^{\infty}(U(1)))$. The line bundle over which we define $A_{i}^{sf}$ is a twisted line bundle twisted by this class. We denote this bundle by $V_{z}^{i}$. Note that here we  have two bundles $V_{z}^{1}, V_{z}^{2}$. There is one subtlety. In fact $V_{z}^{i}$ does not exist globally on $\mathcal{M}$ as a twisted bundle as there are sheet-changing monodromies. What does exists globally over $\mathcal{M}$ is the direct sum of   $V_{z}^{1}$ and $V_{z}^{2}$.

 It seems the following conjectures can be extracted from \citep{GMN1}.
\begin{con} The twisted  bundle $V_{z}$   has a $C^{\infty}$ splitting \begin{equation}V_{z}=V_{z}^{1}\oplus V_{z}^{2}\end{equation}\end{con}

This conjecture brings some importance to $V_{z}^{i}$. But even without the conjecture $V_{z}^{i}$ is still very interesting because of the construction of hyperholomorphic connections on it (see below).

\begin{dfn}The connection $A_{i}^{sf}$ is called the semiflat connection on $V_{z}^{i}$. By taking the direct sum $A_{i}^{sf}$ (for all $i$) induce a connection $A^{sf}$ called the semiflat connection on the direct sum of $V_{z}^{i}$'s.\end{dfn}

Let $\theta_{\gamma}$ be the angular function (coordinate) on the fiber of Hitchin's fibration corresponding to $\gamma$. Write $\gamma_{i}$ as $\gamma_{i}=\gamma_{i}^{0}+\gamma$. Let $\xi$ be the twistor parameter in section 2.2, we define  locally define functions $\mathcal{Y}_{\gamma}^{sf}$ on $\mathcal{M}$ for each $\xi$ by $$\mathcal{Y}_{\gamma}^{sf}:=\exp[\pi RZ_{\gamma}/\xi+i\theta_{\gamma}+\pi R\bar{Z}_{\gamma}\xi]$$ $$\mathcal{Y}_{i}^{sf}:=\exp[\pi RZ_{\gamma_{i}^{0}}/\xi+\pi R\bar{Z}_{\gamma_{i}^{0}}\xi]$$ \begin{equation}\mathcal{Y}_{\gamma_{i}}^{sf}:=\mathcal{Y}_{i}^{sf}\mathcal{Y}_{\gamma}^{sf}\end{equation}
The twistor construction of hyperkahler metrics can be applied to $\mathcal{Y}_{\gamma}^{sf}$ and produce a hyperkahler metric. This is called the $semiflat \ metric$. We will describe the semiflat metric explicitly in section 4 below. For this section we only need to know that the semiflat metric is nonsingular when it is restricted to the open subset $\mathcal{M}^{reg}$ of the Hitchin's moduli space consisting of regular fibers of the Hitchin's fibration. \begin{thm}\citep{GMN1}The semiflat connection $A_{i}^{sf}$ ($A^{sf}$) is hyperholomorphic over $\mathcal{M}^{reg}$ endowed with the semiflat metric.\end{thm}
\begin{thm}\citep{GMN1}For any $\xi$, with respect to the holomorphic structure of $V_{z}^{i}$ induced by the hyperholomorphic connection  $A_{i}^{sf}$, $\mathcal{Y}_{i}^{sf}$  is a  holomorphic section.\end{thm}

We now have $\Gamma_{i}$ and $Z_{\gamma_{i}}$, we still need $\omega(\gamma, \gamma_{i})$. There is no direct definition of them corresponding to the counting of some trajectories. There is a tentative definition of them in section 3 of \citep{GMN1} in terms of coefficients of leading logarithmic divergent terms of $Z_{\gamma_{i}}$ near vanishing cycles.

Suppose we have a set of BPS numbers $\omega(\gamma,\gamma_{i})$ compatible with the BPS numbers $\Omega(\gamma), \omega(\gamma_{ij},\gamma)$. In other words we have $\omega: \Gamma\times\coprod_{i}\Gamma_{i}\rightarrow \mathbf{Q}$ satisfying (2),(3)\footnote{Anticipating the relation with $\mathbf{Z}_{2}$ gerbes one might want to insist that they are all half integers. Anyway the presentation below does not require this.}. Gaiotto-Moore-Neitzke propose a definition of $\mathcal{Y}_{\gamma_{i}}$ by integral equations. First we define a section $x_{\gamma_{i}}(\xi)$ of $V_{z}^{i}$ for each $\xi$ by
\begin{equation}x_{\gamma_{i}}(\xi):=\mathcal{Y}_{\gamma_{i}}^{sf}(\xi)\exp[{1\over 4\pi i}\sum_{\gamma^{'}}\omega(\gamma^{'}, \gamma_{i})\int_{l_{\gamma^{'}}}{d\xi^{'}\over \xi^{'}}{\xi^{'}+\xi\over \xi^{'}-\xi}\log (1-\mathcal{Y}_{\gamma^{'}}(\xi^{'}))]\end{equation}where the sum runs over all charges. Then we define linear maps $g_{i}: V_{z}^{i}\rightarrow V_{z}^{1}\oplus V_{z}^{2}$ ($= V_{z}$ if conjecture 2.3 is true) as a solution of the following equations\footnote{Of course we do not know if there is only one solution with a prescribed asymptotic behavior. Without the uniqueness the definition of $\mathcal{Y}_{\gamma_{i}}$ and hence the construction of hyperholomorphic connection described in this section  is not canonical. This possibility does not mean any troubles. On the other hand there should be a completely geometric definition of $\mathcal{Y}_{\gamma_{i}}$ which is unique. See below for discussions. This would give us a hyperholomorphic connection with more geometric meanings.} \begin{equation}g_{i}(\xi):=g_{i}^{sf}(\xi)+ {1\over 4\pi i}\sum_{\gamma_{li},\ l\neq i}\mu(\gamma_{li}) \int_{l_{\gamma_{li}}}{d\xi^{'}\over \xi^{'}}{\xi^{'}+\xi\over \xi^{'}-\xi}g_{l}(\xi^{'})x_{\gamma_{li}}(\xi^{'})\end{equation}where \begin{equation}x_{\gamma_{ij}}:=\sigma(\gamma_{ij},\gamma_{j})^{-1}x_{\gamma_{i}}x_{\gamma_{j}}^{-1}\end{equation}and $g_{i}^{sf}$ is the inclusion map. Finally we define \begin{equation}\mathcal{Y}_{\gamma_{i}}:=g_{i}x_{\gamma_{i}}\end{equation}

 It is easy to see that  $\mathcal{Y}_{\gamma_{i}}$ is piecewise holomorphic with respect to $\xi$. By Cauchy's theorem  there are nontrivial transformations if we cross the BPS $\mathcal{K}$-rays and $\mathcal{S}$-rays and they are given by $\mathcal{K}$-factors and $\mathcal{S}$-factors. Now let us vary the central charge. We demand that the 2d-4d wall crossing formula is true for $\mathcal{Y}_{\gamma_{i}}$.  Then by the uniqueness of the factorization of a given product in a given stability chamber (see section 2.1) we can use the 2d-4d wall crossing formula to determine all BPS numbers for all stability chambers of the space of central charges. The following theorem then follows  trivially from the definition and theorem 2.3. \begin{thm}Suppose we have a set of BPS numbers $\omega(\gamma,\gamma_{i})$ compatible with the BPS numbers $\Omega(\gamma), \omega(\gamma_{ij},\gamma)$. The data consisting of $\mathcal{V}_{z}$, $\Gamma$, $\{Z_{\gamma}\}$, $\{\mathcal{Y}_{\gamma}\}$, $\{\Omega(\gamma)\}$, $\{\Gamma_{ij}\}$, $\{Z_{\gamma_{ij}}\}$, $\{\omega(\gamma,\gamma_{ij})\}$, $\{\mathcal{Y}_{\gamma_{ij}}\}$, $\{\Gamma_{i}\}$, $\{Z_{\gamma_{i}}\}$,  $\{\omega(\gamma, \gamma_{i})\}$, $\{\mathcal{Y}_{\gamma_{i}}\}$ and twisting functions satisfies the 2d-4d wall crossing formula.\end{thm}

The integral equation formulation also guarantees the existence of the limit of $\mathcal{Y}_{\gamma_{i}}(\mathcal{Y}^{sf}_{\gamma_{i}})^{-1}$ as $\xi$ goes to $\infty$ (prescribing the asymptotic behavior). So we have a Riemann-Hilbert problem formulation of our data.

 The above construction makes sense for any compatible assignment (i.e. (2)(3) hold) $\omega(\gamma,\gamma_{i})$ in the 2d-4d-wall crossing formula. However we have a geometric definition of $\mathcal{Y}_{\gamma_{ij}}$ given before and only an analytic definition of $\mathcal{Y}_{\gamma_{i}}$. Ideally we should be able to find a geometric definition of $\omega(\gamma,\gamma_{i}), \mathcal{Y}_{\gamma_{i}}, x_{\gamma_{i}}, g_{i}$ which also satisfy the integral equations and the 2d-4d wall crossing formula.

Recall that when we use the twistor construction to build the hyperkahler metric we demand that $\mathcal{Y}_{\gamma}$ is holomorphic in the complex structure $J_{\xi}$ for each $\xi$. Now we do a similar thing. We demand that $\mathcal{Y}_{\gamma_{i}}$ is a  holomorphic section of $V_{z}^{1}\oplus V_{z}^{2}$ for all $\xi$ with holomorphic structures determined by a hyperholomorphic connection $A$. From this requirement we can define a hyperholomorphic connection $A$. This is the analog for hyperholomorphic connections of the twistor construction of hyperkahler metrics. In other words, we have the following theorem proved in section 5 and appendix E of \citep{GMN1}:
\begin{thm}$\{\mathcal{Y}_{\gamma_{i}}\}$ determines a hyperholomorphic connection $A$ on $V_{1}^{1}\oplus V_{z}^{2}$ via the twistor construction. \end{thm}

If  conjecture 2.3 is true then we have constructed a hyperholomorphic connection on $v_{z}$. As pointed out by \citep{GMN1} one amazing thing about this construction is that the splitting is clearly not holomorphic with respect to the holomorphic structures induced by the tautological construction of universal bundles. Of course we can also consider the semiflat connections $A^{sf}_{i}$ and take the direct sum of them to get a (semiflat) hyperholomorphic connection $A^{sf}$. The splitting is holomorphic with respect to this semiflat hyperholomorphic connections. If we stick to the study of algebraic (or complex) geometry of the Hitchin's moduli spaces then we can not obtain the information contained in $A$. Only $A^{sf}$ which is essentially a tautological construction can be revealed. Therefore the Gaiotto-Moore-Neitzke construction gives us an opportunity to go beyond complex geometry  in the study of mirror symmetry of Hitchin's moduli spaces and this process exhibits wall crossing phenomena.

\section{Singular Mirror Fibers of A Mirror Pair}
\subsection{SYZ Mirror Symmetry of Hitchin's Moduli Spaces}
It is known that Strominger-Yau-Zaslow (SYZ) type mirror symmetry exists for Hitchin's moduli spaces with dual gauge groups. In general the SYZ mirror conjecture \citep{SYZ,H3}  predicts
\begin{con}Let $X$ be a Calabi-Yau manifold with a B field (denoted by $\mathcal{B}$) and $\check{X}$ be its mirror. Then there should a possibly singular special Lagrangian torus fibration $X\rightarrow B$ for some base $B$.  $\check{X}$ should also be the total space of  a possibly singular special Lagrangian torus fibration over the same base. And this should be a dual fibration in the sense that $\check{X}$ is the moduli space of pairs $(Z,T)$ where $Z$ is a special Lagrangian torus in $X$ and $T$ is a flat trivialization of  $\mathcal{B}$ \end{con}
This version of conjecture is not very precise and probably needs modifications. However for Hitchin's moduli spaces we now have a much better understanding of it for two reasons. First $\mathcal{M}$ is hyperkahler and has the Hitchin's fibration whose fibers are holomorphic Lagrangian submanifolds in complex structure $J_{3}$. We can rotate the complex structure to $J_{1}$ (this is known as a $hyperkahler\ rotation$).  From (19) we see that these torus fibers are now special Lagrangian. Second we have a good understanding of the fibers of the Hitchin's fibration so that we can explicitly discuss the duality of fibers.

Let $G$ be $SL(2)$ whose compact form is $SU(2)$ (sometimes we  call $\mathcal{M}$ the $SL(2)$ Hitchin's moduli space instead of the $SU(2)$ Hitchin's moduli space --- two names for the same object), then its Langlands dual $^{L}G$ is $SO(3)$. Let $\bar{\mathcal{M}}$ be the $SO(3)$ Hitchin's moduli space on $C$. Then we have
\begin{thm}$\mathcal{M}$ and $\bar{\mathcal{M}}$ are SYZ mirror partners.\end{thm}
\begin{dfn}Two Calabi-Yau $n$-orbifolds $M$ and $\bar{M}$ equipped with B fields $\mathcal{B}$ and $\bar{\mathcal{B}}$ respectively, are SYZ mirror partners if  the regular fibers of the two special Lagrangian fibrations $L_{x}$ and $\bar{L}_{x}$ (over the same point $x$) are dual to each other in the sense that\begin{equation}L_{x}=Triv^{U(1)}(\bar{L}_{x},\bar{\mathcal{B}}),\ \ \bar{L}_{x}=Triv^{U(1)}(L_{x},\mathcal{B})\end{equation}Here it is required that the restriction of B field to the regular fiber is isomorphic to the trivial B field and $Triv^{U(1)}(L_{x},\mathcal{B})$ is the space of equivalence classes of the flat trivializations of the restriction of $\mathcal{B}$.\end{dfn}The theorem is proved in \citep{HT} and generalized in \citep{BD} to cover Hitchin's moduli spaces with certain singularities. There are also generalizations for other gauge groups \citep{DP}.  Since in general the SYZ mirror symmetry exchanges a reductive group and its Langlands dual it is  expected that the so-called geometric Langlands duality can be understood in this way.

The precise meaning of the theorem is the following. Both moduli spaces have Hitchin's fibrations over the same base. We view them as moduli space of semistable Higgs bundles. It turns out that when we define the moduli spaces we can prescribe the degree of the bundle part of Higgs bundles. Let us fix this topological type to be degree $d$ and denote the corresponding moduli spaces of semistable Higgs bundles by $\mathcal{M}_{Dol}^{d}$ and $\bar{\mathcal{M}}_{Dol}^{d}$. We denote the corresponding Hitchin fibers by $P^{d}$ and $\bar{P}^{d}$. $\mathcal{M}_{Dol}^{d}$ carries a flat gerbe (B-field) $\mathcal{B}$. In fact as long as $d$ and the rank $r$ of the bundle (in our case $r=2$) are coprime one can show by algebraic geometric methods that the universal Higgs bundle exists over $\mathcal{M}_{Dol}^{d}\times C$. However for the bundle part of the Higgs bundles because of nontrivial automorphism the corresponding universal bundle does not exist. We have discussed this issue in the previous section and we know that what we have is only a twisted universal bundle twisted by a flat gerbe (B-field) $\mathcal{B}$.  We can also define a B-field (see \citep{HT}) $\bar{\mathcal{B}}$ on $\bar{\mathcal{M}}_{Dol}^{d}$. It is proved that the restriction of $\mathcal{B}$ to $P^{d}$ is a trivial gerbe. Let $\mathcal{B}^{e}$ be the $e$-fold tensor product of $\mathcal{B}$ and $Triv^{U(1)}(P^{d},\mathcal{B}^{e})$ the space of equivalence classes of the flat trivializations of the restriction of $\mathcal{B}^{e}$ as a trivial $U(1)$ gerbe. Then there is a smooth identification \begin{equation}Triv^{U(1)}(P^{d},\mathcal{B}^{e})\simeq \bar{P}^{e}\end{equation}And this relation is reciprocal: \begin{equation}Triv^{U(1)}(\bar{P}^{d},\bar{\mathcal{B}}^{e})\simeq P^{e}\end{equation} In this sense we say $\mathcal{M}_{Dol}^{d}$ equipped with $\mathcal{B}^{e}$  and  $\bar{\mathcal{M}}_{Dol}^{e}$ equipped with $\bar{\mathcal{B}}^{d}$ are SYZ mirror partners. This duality is also known as the $fiberwise \ T-duality$.

\subsection{T-duality and Singular Fibers}

The fiberwise T-duality described in the previous section are for regular fibers. What happens for singular fibers is a very important problem. Of course to study this problem one needs to understand the singular fibers and their mirrors first.

The following mirror pairs of singular fibers of Hitchin's fibrations are described in \citep{FW}.

Let $C$ be a smooth elliptic curve whose affine equation is \begin{equation}f(x)=x^{3}+ax+b\end{equation} Let $p$ be the points at the infinity. Let $\mathcal{M}^{0}$ be the $SL(2)$ Hitchin's moduli space with exactly one regular singularity. We assume the degree of the bundle is zero. So the associated moduli space of stable Higgs bundles is $\mathcal{M}_{Dol}^{0}$. Without loss of generality the singularity is $p$. To specify the moduli problem we need to specify the asymptotic behaviors of solutions of Hitchins' equations near $p$. It is described in \citep{FW} but we do not need that information\footnote{But we do assume that the parameters in the specification of asymptotic behaviors are generic. See \citep{KW} for a discussion of this kind of issues.}. Denote the ramified Higgs bundles by $(E,\varphi)$. Then the Hitchin's fibration $\pi$ is given by \begin{equation}\pi((E,\varphi))=Tr\ \varphi^{2}\end{equation}

$\mathcal{M}^{0}$ has complex dimension $2$. As a complex variety (surface) $\mathcal{M}_{Dol}^{0}$ has been determined explicitly \citep{FW}. Its equation is (we just write the affine equation):\begin{equation}\rho^{2}=-(2u+w)f(u)+{f^{'}(u)^{2}\over 4}\end{equation}We refer the readers to section 3 of \citep{FW} for the meaning of $\rho,u,w$. The Hitchin's fibration is given by \begin{equation}(\rho,u,w)\rightarrow w_{0}:=-w\sigma_{0}^{2}/2\end{equation}where $\pm\sigma_{0}$ are eigenvalues of the residue of $\varphi$ at $p$.

Let $F_{w}$ be the fiber at $w$. The fiber is singular if two roots of the right hand side of (62) coincide. It is easy to show that there are three such values of $w$ denoted by $e_{i}, i=1,2,3$. So there are three singular fibers $F_{e_{i}}, i=1,2,3$. Each of them is a curve of arithmetic genus 1 which is a pair of smooth genus 0 curves meeting  at two double points. So each singular fiber has two $A_{1}$ singularities.

There is a group $Q=\mathbf{Z}_{2}\times\mathbf{Z}_{2}$ acting on our moduli space. Let $T_{1}, T_{2}, T_{3}$ be the three nontrivial elements of $Q$. The action of $T_{1}$ is $$u\rightarrow {e_{1}u+e_{2}e_{3}-e_{1}e_{3}-e_{1}e_{2}\over u-e_{1}}$$ $$w\rightarrow w$$\begin{equation}\rho\rightarrow -{(e_{1}-e_{2})(e_{1}-e_{3})\over (u-e_{1})^{2}}\rho\end{equation} $T_{2}$ and $T_{3}$ are given by taking cyclic permutations of $e_{1}, e_{2}, e_{3}$. Let $F_{e_{1}}^{\pm}$ be the two components of the singular fiber $F_{e_{1}}$. $T_{1}$ maps each component $F_{e_{1}}^{\pm}$ to itself and fixes the two double points. $T_{2}, T_{3}$ exchange the two components and the two double points. The action of $Q$ has a moduli interpretation as twisting the bundle part of a Higgs bundle by a line bundle of order 2 while keeping the Higgs field.

We can also consider the $SL(2)$ Hitchin's moduli space with nontrivial $det\ E$. Let $L=det \ E$. As pointed out in \citep{FW}, a twisting of $L$ can be compensated by a twisting of $E$ and hence modulo this action we only care about the degree modulo 2. The corresponding moduli space is denoted by $\mathcal{M}^{1}$ and the associated moduli space of stable Higgs bundles is\footnote{In \citep{FW}, it is denoted by $\mathcal{M}_{H}(C; SL_{2}^{*})$.}$\mathcal{M}^{1}_{Dol}$.

Let $\mathcal{M}$ be the $SL(2)$ Hitchin's moduli space consisting of the two components $\mathcal{M}^{0}$ and $\mathcal{M}^{1}$. According to the SYZ mirror symmetry of Hitchin's moduli spaces in the previous section, the SYZ partner of $\mathcal{M}$ is the $SO(3)$ Hitchin's moduli space of $\bar{\mathcal{M}}$. It consists of two components: $\bar{\mathcal{M}}^{0}$ and $\bar{\mathcal{M}}^{1}$.

Since $SO(3)=SL(2)/\mathbf{Z}_{2}$ is the adjoint form of $SL(2)$ and the three dimensional representation of $SO(3)$ is the adjoint representation, we can understand $SO(3)$ bundles in terms of $SL(2)$ bundle. The moduli space of stable $SO(3)$ Higgs bundles has two components associated to the second Stieffel-Whitney class $w_{2}(W)$ of the underlying bundle $W$. If $w_{2}(W)=0$ then we get $\bar{\mathcal{M}}^{0}$. If $w_{2}(W)$ is the nontrivial element of $H^{2}(C,\mathbf{Z}_{2})=\mathbf{Z}_{2}$ then we get $\bar{\mathcal{M}}^{1}$.

If $w_{2}(W)=0$ then the structure group can be lifted to $SL(2)$ and $W$ can be identified with the adjoint bundle of an $SL(2)$ bundle $E$ determined up to twisting by a line bundle by order 2. So we have \begin{equation} \bar{\mathcal{M}}^{0}=\mathcal{M}^{0}/Q\end{equation}

The action $Q$ is the same one described before. So the three singular fibers of $\mathcal{M}^{0}$ descents to three singular fibers of $\bar{\mathcal{M}}^{0}$. $T_{2}$ (or $T_{3}$) identifies the two components $F^{\pm}_{e_{1}}$ giving us a genus zero curve with two points identified (so we have one double point). The $T_{1}$ quotient gives us
the same genus zero curve with one double point (one $A_{1}$ singularity). This is the description of the corresponding singular fiber of $\bar{\mathcal{M}}^{0}$.

There is a similar relation between $\mathcal{M}^{1}$ and $\bar{\mathcal{M}}^{1}$\begin{equation} \bar{\mathcal{M}}^{1}=\mathcal{M}^{1}/Q\end{equation} There is also a careful analysis of B-fields in in terms of spin structures on $C$ and square roots of line bundles over the moduli spaces in \citep{FW}.
We know that the parts of $\mathcal{M}$ and $\bar{\mathcal{M}}$ consisting of regular fibers satisfy fiberwise T-duality. The above analysis also gives us a canonical way to associate a singular fiber of $\mathcal{M}$ (as a holomorphic Lagrangian fibration or a special Lagrangian fibration after a hyperkahler rotation) to a singular fiber of $\bar{\mathcal{M}}$ and vice versa. So even though in general it is not very clear what the mirror of a Calabi-Yau orbifold which is the total space of a singular special Lagrangian torus fibration should be exactly (because we do not know what the fiberwise T-duality of a general singular torus fiber means), at least for Hitchin's moduli spaces in this paper we should take the full SYZ mirror of $\mathcal{M}$ to be $\bar{\mathcal{M}}$ and vice versa (with B-fields). In this sense we identify the singular fibers of $\bar{\mathcal{M}}$ described above as the SYZ mirrors of singular fibers of $\mathcal{M}$ and vice versa.

It is then important to understand what the singular fiberwise T-duality is for these mirror pairs of singular fibers of mirror pairs of Calabi-Yau orbifolds (our Hitchin's moduli spaces). In general regular fiberwise T-duality not only identifies fibers as spaces of trivializations of B-fields but also transforms skyscraper sheaves (B-branes) supported as a point in a fiber to Lagrangian A-brane consisting of the mirror fiber endowed with a flat unitary line bundle (if we ignore the B-fields). Frenkel and Witten studied the extension of this picture to singular mirror fibers described above. They have also understood the B-fields's effects on singular fibers. The singular T-duality analyzed in their paper is interpreted as a so-called endoscopy phenomenon. This is a very deep aspect of Langlands duality.

Another thing that a regular fiberwise T-duality does is to transform certain connections (see section 5) to special Lagrangian sections endowed with appropriate connections. This is more complicated and one of the  purposes of this paper is to prepare for a study of the extension of this phenomenon to singular mirror fibers.

\section{Ooguri-Vafa Spaces}
\subsection{Ooguri-Vafa Metrics}

We describe in this section some nice hyperkahler spaces that will be used later. We follow closely the exposition of \citep{GMN1}.
Let $D_{\Lambda}^{*}$ be a punctured disk of the complex plane \begin{equation}D_{\Lambda}^{*}:=\{a\in \mathbf{C}\mid0<|a|<|\Lambda|\}\end{equation}where $\Lambda$ is a constant. Let $D_{\Lambda}$ be the associated disk with the center $a=0$. We define a rank 2 local system of lattices $\Gamma$ with an integral antisymmetric pairing over $D_{\Lambda}^{*}$ by going to its universal cover $C_{\Lambda}:=\{z\mid Re(z)<\log |\Lambda|\}, e^{z}=a$ and setting the pull-back of the local system (also denoted by $\Gamma$) as\begin{equation}\Gamma:=(C_{\Lambda}\times (\mathbf{Z}_{\gamma_{e}}\oplus\mathbf{Z}_{\gamma_{m}}))/\mathbf{Z}\end{equation} where $\gamma_{e}, \gamma_{m}$ are the two generators. The pairing is defined by linearity and \begin{equation}\langle\gamma_{e}, \gamma_{m}\rangle=1\end{equation}The generator $1\in\mathbf{Z}$ acts by (denote the action by $T$)\begin{equation}T((z,q\gamma_{e}+p\gamma_{m})):=(z+2\pi i, (q-\Delta p)\gamma_{e}+p\gamma_{m})\end{equation}where $q,p,\Delta$ are integers. Clearly $\Delta$ determines the monodromy of $\gamma_{m}$ and $\gamma_{e}$ has trivial monodromy.

The central charges are defined by linearity and \begin{equation}Z(\gamma_{e}):=a=e^{z}\end{equation}\begin{equation}Z(\gamma_{m}):=\tau_{0}a+{\Delta\over2\pi i}(a\log{a\over \Lambda}-a)\end{equation}where $\tau_{0}$ is a constant with positive imaginary part.  When $\Delta\neq0$, we can redefine $\Lambda$ to absorb $\tau_{0}$. So when $\Delta\neq0$ we can write \begin{equation}Z(\gamma_{m}):={\Delta\over2\pi i}(a\log{a\over \Lambda}-a)\end{equation}$Z_{\gamma_{e}}$ and $Z_{\gamma_{e}}$ will play the roles of $a$ and $a_{D}$ in section 2.2. Therefore we define \begin{equation}\tau(a):= {\Delta\over 2\pi i}\log{a\over \Lambda}\end{equation}

We consider the dual local system $\Gamma^{*}\otimes\mathbf{R}/2\pi \mathbf{Z}$ of the local system $\Gamma$ over the universal cover $C_{\Lambda}$\begin{equation}\Gamma^{*}\otimes\mathbf{R}/2\pi \mathbf{Z}:=N_{\Lambda}/\mathbf{Z}\end{equation}\begin{equation}N_{\Lambda}:=C_{\Lambda}\times ((\mathbf{R}/2\pi \mathbf{Z})\gamma_{e}^{*}\oplus(\mathbf{R}/2\pi \mathbf{Z})\gamma_{m}^{*})\end{equation}where the $\mathbf{Z}$ action is \begin{equation}T((z, \theta_{e}\gamma_{e}^{*}+\theta_{m}\gamma_{m}^{*})):=(z+2\pi i,\theta_{e}\gamma_{e}^{*}+(\theta_{m}+\Delta\theta_{e})\gamma_{m}^{*})\end{equation} We use $\mathcal{M}^{sf}$ to denote $\Gamma^{*}\otimes\mathbf{R}/2\pi \mathbf{Z}$. It is the total space of a fibration over $D_{\Lambda}^{*}$. $\theta_{e}, \theta_{m}$ are coordinates of the fibers. Each fiber is a torus.

For later convenience we follow \citep{GMN1} and modify our representation a little bit by setting the $\mathbf{Z}$ action as \begin{equation}\hat{T}((z, \theta_{e}\gamma_{e}^{*}+\theta_{m}\gamma_{m}^{*})):=(z+2\pi i,\theta_{e}\gamma_{e}^{*}+(\theta_{m}+\Delta\theta_{e}+\pi\Delta)\gamma_{m}^{*})\end{equation}and define\footnote{The reason of doing this shift is explained in \citep{GMN3} and has something to do with the quadratic refinements.} \begin{equation}\hat{\mathcal{M}}^{sf}:=N_{\Lambda}/\mathbf{Z}\end{equation}where $\mathbf{Z}$ acts by $\hat{T}$. The two spaces $\hat{\mathcal{M}}^{sf}$ and $\mathcal{M}^{sf}$ are identified by a shift of $\theta_{e}$.

$\hat{\mathcal{M}}^{sf}$ can carry the following hyperkahler  metric called the semiflat metric\begin{equation}g^{sf}:=R Im \tau|da|^{2}+{1\over 4\pi^{2}R Im\tau}|d\theta_{m}-\tau d\theta_{e}|^{2}\end{equation}where $R$ is a real constant. The definition is invariant under $\hat{T}$ and hence descends to $\hat{\mathcal{M}}^{sf}$. It is positive definite provided $\Delta>0$. The restriction to each fiber is a standard flat metric on the torus. The base part of the metric cannot be extended to $D_{\Lambda}$ as $\tau$ is singular at $a=0$.

On the other hand we can also view $\hat{\mathcal{M}}^{sf}$ as a principal $U(1)$ bundle over $D_{\Lambda}^{*}\times S^{1}$ by treating $\theta_{m}$ as the fiber coordinate. This can happen because $\theta_{e}$ is invariant under $\hat{T}$. The first Chern class of this $U(1)$ bundle is $-\Delta$.

The semiflat metric is not the only possible hyperkahler metric on the total space of this $U(1)$ bundle. More of them can be obtained from the Gibbons-Hawking Ansatz. For a principal $U(1)$ bundle over a region of $\mathbf{R}^{3}$ this ansatz is \begin{equation}U^{-1}({d\chi\over2\pi}+A_{GH})^{2}+Ud\vec{x}^{2}\end{equation}where $d\vec{x}:=(x_{1}.x_{2},x_{3})\in \mathbf{R}^{3}$ and $\chi$ is the coordinate of the fiber. $U$ is a harmonic function and \begin{equation}dA_{GH}=\star dU\end{equation}where $\star$ is the Hodge dual computed with respect to the standard flat metric on $\mathbf{R}^{3}$.

To adapt it to our situation we let the third coordinate $x_{3}$ of the Gibbons-Hawing ansatz to be periodic and set \begin{equation}(x_{1}+ix_{2},x_{3})=(a, {\theta_{e}\over 2\pi R}),\ \ \chi=\theta_{m}\end{equation}The semiflat metric fits into this ansatz by setting \begin{equation}U^{sf}=R Im\tau,\ \ A_{GH}^{sf}=-R Re\tau dx_{3}\end{equation}

The Ooguri-Vafa metric is obtained by setting \begin{equation}U=U^{\Omega}:=\sum_{q=1}^{\infty}q^{2}\Omega_{q}U_{q}\end{equation}Here $\Omega_{q}$ is an integer and \begin{equation}\Omega_{q}=\Omega_{-q}\end{equation}$\Delta$ is determined by them via\begin{equation}\Delta={1\over 2}\sum_{q}q^{2}\Omega_{q}=\sum_{q>0}q^{2}\Omega_{q}\end{equation}$U_{q}$ is the following harmonic function \begin{equation}U_{q}:={1\over 4\pi}\sum_{n=-\infty}^{\infty}({1\over\sqrt{q^{2}|a|^{2}+R^{-2}(q\theta_{e}/2\pi+n)^{2}}}-\kappa_{n})\end{equation}where $\kappa_{n}$ is a constant determined by the condition that leading term in $U$ for large $R$ is $U^{sf}$. By Poisson summation $U_{q}$ can be rewritten (see \citep{GW0}) as \begin{equation}U_{q}=U_{q}^{sf}+U^{inst}_{q}\end{equation}\begin{equation}U_{q}^{sf}=-{R\over 4\pi}(\log{a\over \Lambda}+\log{\bar{a}\over\bar{\Lambda}})\end{equation}\begin{equation}U^{inst}_{q}={R\over 4\pi}\sum_{n\neq0}e^{inq\theta_{e}}K_{0}(2\pi R|nqa|)\end{equation}where $K_{0}$ is the modified Bessel function. Then it is not hard to show that there is a constant $C$ such that \begin{equation}|U_{q}-U^{sf}_{q}|\leq CR e^{-2\pi R|qa|}\end{equation} Note that $U^{sf}$ is obtained by summing $U^{sf}_{q}$. Therefore the deviation of the Ooguri-Vafa metric from the semiflat metric is exponentially suppressed for large $R$. In fact the Ooguri-Vafa metric can be defined on the total space of the same principal $U(1)$ bundle with Chern class $-\Delta$ (equation (87) guarantees this) defined before over $D_{\Lambda}^{*}\times S^{1}$. When $R\rightarrow \infty$, it goes to the semiflat metric.

Once we have $U^{\Omega}$, the one form $A_{GH}$ is determined by (82). We denote it by\begin{equation}A^{\Omega}(=A_{GH}) =A^{\Omega}_{a}da+A^{\Omega}_{\bar{a}}d\bar{a}+A^{\Omega}_{\theta_{e}}d\theta_{e}\end{equation}It has no $\theta_{m}$ component. More explicitly,\begin{equation}A^{\Omega}=A^{\Omega}_{sf}+A^{\Omega}_{inst}\end{equation}\begin{equation}A^{\Omega}_{sf}=A^{\Omega}_{\theta_{e}}d\theta_{e}=
\sum_{q=1}^{\infty}q^{2}\Omega_{q}A_{q}^{sf}\end{equation}\begin{equation}A_{q}^{sf}={i\over 8\pi^{2}}(\log{a\over\Lambda}-\log{\bar{a}\over\bar{\Lambda}})d\theta_{e}\end{equation}\begin{equation} A^{\Omega}_{inst}=A^{\Omega}_{a}da+A^{\Omega}_{\bar{a}}d\bar{a}=\sum_{q=1}^{\infty}q^{2}\Omega_{q}A_{q}^{inst}\end{equation}
\begin{equation}A_{q}^{inst}=-{R\over 4\pi}({da\over a}-{d\bar{a}\over \bar{a}})\sum_{n\neq 0}sgn(n)e^{inq\theta_{e}}|a|K_{1}(2\pi R|nqa|)\end{equation}

The semiflat metric  can not be extended to $a=0$ but for  the Ooguri-Vafa one can do better. In fact define \begin{equation}s_{n,q}:=(a=0, \theta_{e}=2\pi n/q),\ \ S_{q}:=\{s_{n,q}\mid n\in\mathbf{Z}\}\end{equation} then $U^{\Omega}$ is well-defined outside $S_{\Omega}:=\cup_{q:\Omega_{q}\neq 0}S_{q}$. Denote the total space of the $U(1)$ bundle over $D_{\Lambda}^{*}\times S^{1}$ by $\mathcal{M}(\Omega)$ (for the semiflat metric we used the notation $\hat{\mathcal{M}}^{sf}$. They are the same manifold with  different metrics). Then we have the projection\begin{equation}\pi_{\Omega}:\mathcal{M}(\Omega)\rightarrow D_{\Lambda}^{*}\times S^{1}\end{equation}which clearly can be extended to \begin{equation}\pi_{\Omega}:\mathcal{M}(\Omega)\rightarrow D_{\Lambda}\times S^{1}-S_{\Omega}\end{equation} Here we have used  the same notations for the extension. We fill in all the points $s_{n,q}$ to $D_{\Lambda}\times S^{1}-S_{\Omega}$ to get $D_{\Lambda}\times S^{1}$ and do the extension in the fiber direction by adding a single point $p_{n,q}$ to $\mathcal{M}(\Omega)$ over each $s_{n,q}$ and demand continuity. Therefore the circle fibers of the $U(1)$ bundle collapses over $s_{n.q}$ to a single point. Denote this extension by $\hat{\mathcal{M}}(\Omega)$. Then $\hat{\mathcal{M}}(\Omega)$ is a singular $U(1)$ bundle over $D_{\Lambda}\times S^{1}$ with finitely many collapsing fibers. Outside these fibers the Ooguri-Vafa metric is nonsingular. More careful analysis shows that locally near each singularity $p_{n,q}$ the metric is asymptotic to an asymptotically locally Euclidean hyperkahler metric for a Kleinian singularity $\mathbf{C}^{2}/\mathbf{Z}_{N_{q}}$ for some integer $N_{q}$ (see section 6 of \citep{GMN1}).

\subsection{Twistor Coordinates and 4d Wall Crossing}

A key observation due to Gaiotto, Moore and Neitzke is that the Ooguri-Vafa metric described above has another description in terms of twistor coordinates such that there is a manifest 4d wall crossing phenomenon involved.

Let us start with the semiflat metric. Define the following so-called twistor coordinates \begin{equation}\mathcal{X}_{e}^{sf}:=\exp[\pi RZ_{{e}}/\xi+i\theta_{e}+\pi R\bar{Z}_{e}\xi]\end{equation}\begin{equation}\mathcal{X}_{m}^{sf}:=\exp[\pi RZ_{{m}}/\xi+i\theta_{e}+\pi R\bar{Z}_{m}\xi]\end{equation}where $\xi\in \mathbf{C}^{\times}$. Define a two form over $\hat{\mathcal{M}}^{sf}$ for every $\xi$ by \begin{equation}\mathbf{\Omega}(\xi):={1\over 8\pi^{2}R}{d\mathcal{X}^{sf}_{e}\over \mathcal{X}^{sf}_{e}}\wedge{d\mathcal{X}^{sf}_{m}\over \mathcal{X}^{sf}_{m}}\end{equation} The twistor  construction  essentially says that $\xi$ is the twistor parameter and $\mathbf{\Omega}(\xi)$ for each $\xi$ is the holomorphic symplectic form of the hyperkahler metric in the complex structure corresponding to $\xi$. From this requirement one can  reconstruct the hyperkahler metric on $\hat{\mathcal{M}}^{sf}$ which is exactly the semiflat metric.

Now we can do the similar thing to the Ooguri-Vafa metric. For the Ooguri-Vafa metric over $\mathcal{M}(\Omega)$ the twistor coordinates are \begin{equation}\mathcal{X}_{e}:=\mathcal{X}_{e}^{sf}\end{equation}\begin{equation}\mathcal{X}_{m}:=\mathcal{X}_{m}^{sf}\prod_{q\in\mathbf{Z},q\neq 0}F_{q}^{q\Omega_{q}}\end{equation}To define $F_{q}$ we go to the universal cover $N_{\Lambda}$ times the universal cover of $\mathbf{C}^{\times}$. Define \begin{equation}F_{q}(z,\theta_{e},\xi):=\exp[-{1\over 4\pi i}\int_{l_{q\gamma_{e}}}{d\xi^{'}\over \xi^{'}}{\xi^{'}+\xi\over \xi^{'}-\xi}\log [1-\mathcal{X}_{e}(\xi^{'})^{q}]]\end{equation}and \begin{equation}l_{q\gamma_{e}}:=\{(z,\xi)\mid {qZ(\gamma_{e})\over\xi}\in \mathbf{R}_{-}\}\end{equation}So  for each fixed $z$, $l_{q\gamma_{e}}$ is a ray. $F_{q}$ is piecewise holomorphic with respect to $\xi$ away from these rays. If one crosses a ray then the above definition tells us that  one needs to make a nontrivial transformation given by the following identity \begin{equation}(\mathcal{X}_{m})_{l_{q\gamma_{e}}}^{+}=(\mathcal{X}_{m})_{l_{q\gamma_{e}}}^{-}\prod (1-\mathcal{X}_{e}^{q})^{q\Omega_{q}}\end{equation}where $(\mathcal{X}_{m})_{l_{q\gamma_{e}}}^{+}$ and $(\mathcal{X}_{m})_{l_{q\gamma_{e}}}^{-}$ are the limits of $\mathcal{X}_{m}$ approaching $l_{q\gamma_{e}}$ from the clockwise and counterclockwise directions respectively. Note that the order of the product does not matter. When we change $a$  the rays will rotate but the order of the product still does not matter. Also note that the space of central charges is the space of $a$.

To handle $F_{q}$'s dependence on $z$ we can analytically continue in $z$ to define a function (also denoted by $F_{q}$) with no discontinuity in $z$. Then $F_{q}$ is quasi-periodic \begin{equation}F_{q}(z+2\pi i,\theta_{e},\xi)=(1-\mathcal{X}^{q}_{e})^{-1}F_{q}(z,\theta_{e},\xi)\end{equation}We can check that it descends and then extends to $\hat{\mathcal{M}}(\Omega)$.

Now we define \begin{equation}\Omega(q\gamma_{e}+p\gamma_{m}):=\delta_{p,0}\Omega_{q}\end{equation}Compare the data consisting of $\mathcal{X}$, $\Gamma$, $Z$, $\Omega(q\gamma_{e}+p\gamma_{m})$ (see \citep{GMN3} for the discussion about twisting functions\footnote{\citep{GMN3} uses the language of quadratic refinements. One needs to translate it to twisting functions by the recipe of section 7 of \citep{GMN1}.}) and the 4d wall crossing formalism in section 2.1 and we get \begin{thm}The data for the twistor description of the Ooguri-Vafa metric given above satisfies the 4d wall crossing formalism. The nontrivial transformations across $l_{q\gamma_{e}}$ (identified as BPS $\mathcal{K}$-rays) are $\mathcal{K}$-factors and the 4d wall crossing formula they satisfy is the trivial one (no change of order). \end{thm}

\subsection{Hyperholomorphic Connections over Ooguri-Vafa Spaces}

In the setting of the Ooguri-Vafa space, we define \begin{equation}W:={\delta\over2\pi i}(a\log{a\over\Lambda}-a)+W^{analytic}\end{equation}where $\delta$ is an integer or half integer and $W^{analytic}$ is an analytic function of $a$: $W=w_{1}a+w_{2}a^{2}+\cdots$. Setting $\Lambda_{1}$ by $w_{1}={\delta\over2\pi i}(\log{\Lambda\over\Lambda_{1}})$ we can write\begin{equation}W={\delta\over2\pi i}(a\log{a\over\Lambda_{1}}-a)+W^{analytic}\end{equation}where \begin{equation}W^{analytic}=w_{2}a^{2}+\cdots\end{equation}Define
\begin{equation}\omega(q\gamma_{e}+p\gamma_{m},\gamma_{0}):=\delta_{p,0}\omega_{q}\end{equation}where $\gamma_{0}$ is an extra formal variable whose meaning will be clarified at the end of this section. $\omega_{q}$ is a rational number and \begin{equation}\omega_{q}=\omega_{-q}\end{equation}They are required to satisfy \begin{equation}\delta={1\over2}\sum_{q\in\mathbf{Z}}q\omega_{q}=\sum_{q>0}q\omega_{q}\end{equation}Define\begin{equation}\eta+\tau\alpha=\partial W/\partial a=t\end{equation}where $\eta,\alpha$ are real. Denote $W-W^{analytic}$ by $W_{0}$. Clearly \begin{equation}t={\delta\over\Delta}\tau+t^{analytic}\end{equation}where $t^{analytic}=\partial W^{analytic}/\partial a$. The one form
\begin{equation}A^{sf}=i(\eta d\theta_{e}+\alpha d\theta_{m})\end{equation}defines a connection on the associated line bundle $V$ of the following principal $U(1)$ bundle (the fiber coordinate is $\psi$) over $\hat{\mathcal{M}}^{sf}$\begin{equation}\hat{P}^{sf}:=(N_{\Lambda}\times U(1))/\mathbf{Z}\end{equation} where the generator of $\mathbf{Z}$ acts by \begin{equation}\hat{T}(z,\theta_{e}\gamma_{e}^{*}+\theta_{m}\gamma_{m}^{*},e^{i\psi})=(z+2\pi i,\theta_{e}\gamma_{e}^{*}+(\theta_{m}+\Delta\theta_{e}+\pi\Delta)\gamma_{m}^{*},e^{i\psi}e^{-i\delta\theta_{e}}(-1)^{n_{\omega}})\end{equation}
where $n_{\omega}=\sum_{q>0}\omega_{q}$ which induces a shift of the fiber coordinates analogous to the shift of $\theta_{m}$ in section 4.1. The principal connection inducing $A^{sf}$ is $\Theta^{sf}=i(d\psi+\eta d\theta_{e}+\alpha d\theta_{m})$.
$A^{sf}$ is  the semiflat connection we discussed before. Its curvature is not zero but nonzero contributions all come from $t^{analytic}$. If $W=W_{0}$ then $A^{sf}=i {\delta\over\Delta}d\theta_{m}$ which is flat.

Now define \begin{equation}\mathcal{X}_{W}^{sf}:=\exp({\pi R\over\xi}W-i\psi+\pi R \xi \bar{W})\end{equation}\begin{equation}\mathcal{X}_{W}:=\mathcal{X}_{W}^{sf}\prod_{q\in\mathbf{Z},q\neq 0}F_{q}^{\omega_{q}}\end{equation}It is shown that $\mathcal{X}_{W}$ descends to a $U(1)$ equivariant function from $(N_{\Lambda}\times U(1)/\mathbf{Z}$ to $\mathbf{C}$. Therefore it is a section of the associated line bundle $V$ over $\hat{\mathcal{M}}^{sf}$. In fact $V$ may be just a twisted bundle twisted by a $\mathbf{Z}_{2}$ B field as (122) involves a choice of square root when $\delta$ is a half integer (section 2.6). Recall that $\hat{\mathcal{M}}^{sf}$ supports not only the semiflat metric but also the Ooguri-Vafa metric. Now we use the Ooguri-Vafa metric. Following the philosophy of section 2.6 we define a hyperholomorphic connection $A$ on $V$ by demanding that the section $\mathcal{X}_{W}$ is holomorphic for all compatible complex structures. In other words $(d+A)\mathcal{X}_{W}$ must be type $(1,0)$ for all compatible complex structures compatible to the Ooguri-Vafa metric.

\begin{dfn}The connection defined in the previous paragraph is called the GMN connection.\end{dfn}
The explicit form of the GMN connection can be found in \citep{GMN1}. The Ooguri-Vafa metric is actually nonsingular on  $\hat{\mathcal{M}}(\Omega)$ minus some isolated points (section 4.1). Similarly the GMN connection can be viewed as a nonsingular connection on $V$ over $\hat{\mathcal{M}}(\Omega)$ minus some isolated points. Details are again in \citep{GMN1}.

We can separate the contributions of the leading terms of $W$ and $W^{analytic}$ (as explained in section 6 it seems to be necessary to do this if we want to use the Ooguri-Vafa space as a local model of some Hitchin's moduli spaces). So we introduce $\Omega^{eff}_{q}:=\omega_{q}/q$ and write \begin{equation}\mathcal{X}^{analytic}:=\exp({\pi R\over\xi}W^{analytic}+\pi R \xi \bar{W}^{analytic})\end{equation}Then we have\begin{equation}\mathcal{X}_{W}=\mathcal{X}_{W}^{\omega}\mathcal{X}^{analytic}\end{equation}where $\mathcal{X}_{W}^{\omega}$ is defined to be $\mathcal{X}_{m}$ in section 4.2 computed with $\Omega_{q}^{eff}$. Clearly all discontinuous jumps are associated to $\mathcal{X}^{\omega}_{W}$ and hence are determined by the two leading terms of $W$. One subtle point is that in the definition of $\mathcal{X}_{m}$ there is a term $i\theta_{m}$ in the exponential. Here we identify $i\theta_{m}$ with $-i\psi$.

The analogy of the exponentially suppressed deviation of the Ooguri-Vafa metric from the semiflat metric is the exponentially suppressed (for large $R$) deviation of the GMN connection from the semiflat connection away from  neighborhoods of singular points of the GMN connection.

 A comparison of this section with section 2.6 tells us that the 2d-4d wall crossing formalism is working here.\begin{thm}\citep{GMN1} The above data defining the GMN connection over the Ooguri-Vafa space satisfies the 2d-4d wall crossing formula. In particular $\mathcal{V}$ is a set of a single element (so we only have $\Gamma_{1}$). An element of the torsor $\Gamma_{1}$ is written as $\gamma_{0}+\Gamma$. There is no $\Gamma_{ij}$ and $\mu(\gamma_{ij})=0$. $Z_{1}$ is nothing but our $W$(the value $Z_{\gamma_{0}}$ is unimportant). BPS numbers are $\Omega$ and $\omega$.\end{thm}

\section{SYZ Mirror Symmetry of Ooguri-Vafa Spaces}
\subsection{Holomorphic Coordinates}
For the hyperkahler metric of the Gibbons-Hawking ansatz, the triple of Kahler forms are given by
\begin{equation}\omega_{1}=dx_{1}\wedge({d\theta_{m}\over 2\pi}+A)+Udx_{2}\wedge dx_{3}\end{equation}
\begin{equation}\omega_{2}=dx_{2}\wedge({d\theta_{m}\over 2\pi}+A)+Udx_{3}\wedge dx_{1}\end{equation}
\begin{equation}\omega_{3}=dx_{3}\wedge({d\theta_{m}\over 2\pi}+A)+Udx_{1}\wedge dx_{2}\end{equation}
The holomorphic symplectic form for $J_{3}$ is \begin{equation}\Omega_{3}=da\wedge ({d\theta_{m}\over 2\pi}+A-{iU\over 2\pi R}d\theta_{e})\end{equation}It turns out that for $J_{3}$ we have a holomorphic fibration which can be locally identified with the holomorphic cotangent bundle $T^{*}B$ of the base. Let us quote the following theorems from \citep{GW1}
\begin{thm}Let $a=x_{1}+ix_{2}$ and $B=D_{\Lambda}$. Denote the total space of the Gibbons-Hawking ansatz by $X$. Then $X$ in the complex structure $J_{3}$ is a holomorphic elliptic fibration over the base $B$. There exist holomorphic coordinates $a, z$, where $a$ is the coordinate of $B$ defined before and $z$ is the holomorphic coordinate along the fiber direction, such that $\Omega_{3}=da\wedge dz$. The definition of $z$ depends on the choice of local section of the fibration (in fact a local section provides base points of the fibers for the integration of the restriction of ${d\theta_{m}\over 2\pi}+A-{iU\over 2\pi R}d\theta_{e}$).\end{thm}
\begin{thm}For the Ooguri-Vafa space the holomorphic elliptic fibration $\hat{\mathcal{M}}(\Omega)\rightarrow B$ has periods $1$ and $\tau$ where $\tau$ is defined in (74). \end{thm}
\subsection{Semiflat SYZ Mirror Symmetry}
Let us consider the semiflat hyperkahler space $\hat{\mathcal{M}}^{sf}$ with the semiflat metric . There is the nonsingular torus fibration \begin{equation}\hat{\mathcal{M}}^{sf}=\Gamma^{*}\otimes\mathbf{R}/(2\pi\mathbf{Z})\rightarrow B\end{equation}The SYZ mirror symmetry in this case is quite simple. On the base $B$ the restriction of the semiflat metric is not changed. But we replace each torus fiber $T$ (note that they are all nonsingular) by its dual torus $\check{T}= Hom (\Gamma^{*}, U(1))$. Denote the space obtained in this way by $\check{\mathcal{M}}^{sf}$. So \begin{equation}\check{\mathcal{M}}^{sf}=\Gamma\otimes\mathbf{R}/(2\pi\mathbf{Z})\rightarrow B\end{equation}
Since each fiber of it is dual to the corresponding fiber of $\hat{\mathcal{M}}^{sf}$, we say $\hat{\mathcal{M}}^{sf}$ and $\check{\mathcal{M}}^{sf}$ are related by the fiberwise T duality.

We still need to specify the hyperkahler metric on $\check{\mathcal{M}}^{sf}$ which is supposed to be dual to the semiflat one on $\hat{\mathcal{M}}^{sf}$. We first rescale the metric by $4\pi^{2}R$ to eliminate the $R$ dependence of the fiber part metric. The rescaled fiber metric of $\hat{\mathcal{M}}^{sf}$ is \begin{equation}4\pi^{2}Rg^{sf}_{fiber}={1\over Im \tau}|dz|^{2}={1\over Im \tau}(d\theta_{m}-\tau d\theta_{e})(d\theta_{m}-\bar{\tau} d\theta_{e})\end{equation}This fiber metric is flat. Now we take its dual flat metric on the dual torus. By this we mean we replace $\theta_{e}, \theta_{m}$ by their dual coordinates $\check{\theta}_{e},\check{ \theta}_{m}$ and then replace the coefficient matrix of $g^{sf}_{fiber}$ by its inverse. So in the end we get the following dual semiflat metric on $\check{\mathcal{M}}^{sf}$:\begin{equation}4\pi^{2}R\check{g}^{sf}=4\pi^{2}R^{2}Im \tau |da|^{2}+{1\over Im \tau}(d\check{\theta}_{e}+\tau d\check{\theta}_{m})(d\check{\theta}_{e}+\bar{\tau} d\breve{\theta}_{m})\end{equation} $\check{\mathcal{M}}^{sf}$ endowed with this semiflat metric $\check{g}^{sf}$ is called the semiflat SYZ mirror of $\hat{\mathcal{M}}^{sf}$.

Clearly $\check{\mathcal{M}}^{sf}$ endowed with $\check{g}^{sf}$ is isometric to  $\hat{\mathcal{M}}^{sf}$ endowed with $g^{sf}$ if we make the identification\begin{equation}\check{\theta}_{m}\rightarrow \theta_{e}, \check{\theta}_{e}\rightarrow -\theta_{m}\end{equation}In this sense we say that $\hat{\mathcal{M}}^{sf}$ is its own mirror.

Now let us consider the transformations of hyperholomorphic connections. By general expectation the SYZ mirror of a hyperholomorphic connection over $\hat{\mathcal{M}}^{sf}$ should be a holomorphic Lagrangian section of $\check{\mathcal{M}}^{sf}\rightarrow B$ in the complex structure $J_{3}$ (i.e. the complex structure in which the dual fibration is holomorphic) together with a flat $U(1)$ connection on it. Here by holomorphic Lagrangian we mean that the restriction of the holomorphic symplectic form $\Omega_{3}$ is zero.  The reason for this expectation is that by the hyperkahler rotation argument a holomorphic Lagrangian section becomes a special Lagrangian section in $J_{1}, J_{2}$ and  then the general philosophy of mirror symmetry predicts this correspondence (the correspondence of B branes and A branes).

This correspondence between hyperholomorphic connections of $\hat{\mathcal{M}}^{sf}$ and holomorphic Lagrangian sections of $\check{\mathcal{M}}^{sf}$ plus flat $U(1)$ connections can be verified for the semiflat connection. We consider the trivial line bundle over $\hat{\mathcal{M}}^{sf}$. We know that $A^{sf}$ is a hyperholomorphic connection over $\hat{\mathcal{M}}^{sf}$. $A^{sf}$ induces canonically a section $(\alpha \theta_{m}+\eta \theta_{e})$ of  $\check{\mathcal{M}}^{sf}\rightarrow B$ (i.e. $\check{\theta}_{e}=\eta,\ \check{\theta}_{m}=\alpha$). Since $W$ is holomorphic it is easy to show that \begin{thm}\citep{GMN1}The section $(\alpha \theta_{m}+\eta \theta_{e})$ is a holomorphic lagrangian section of $\check{\mathcal{M}}^{sf}\rightarrow B$.\end{thm}Therefore we can canonically produce a holomorphic Lagrangian section (of $\check{\mathcal{M}}^{sf}\rightarrow B$) plus the trivial flat $U(1)$ bundle on it from a semiflat hyperholomorphic connection over $\hat{\mathcal{M}}^{sf}$. It is also easy to invert this construction. More complete studies of semiflat mirror symmetry can be found in \citep{A,L1,LYZ}.

\subsection{SYZ Mirror Symmetry of Ooguri-Vafa Spaces}
We want to find the SYZ mirror of the Ooguri-Vafa space $\hat{\mathcal{M}}(\Omega)$. This is actually not a well-defined problem as the SYZ mirror conjecture does not exactly say what happens for singular fibers. However one can choose properties that are true in the semiflat cases and try to generalize them to the Ooguri-Vafa space case.

We want to make sure that the nonsingular torus fibers of the Ooguri-Vafa space and its mirror are dual to each other. One natural way to achieve that is replacing $\theta_{e},\theta_{m}$ by $\check{\theta}_{m}, -\check{\theta}_{e}$. The result is also an Ooguri-Vafa space.  This would guarantee that when we take the semiflat approximation of the original Ooguri-Vafa space and its mirror we can recover the mirror duality of section 5.2. Another strong support of this recipe is from \citep{Ch} where it is proved that the counting of holomorphic disk instantons and the reconstruction procedure of SYZ mirrors produce the same result.

We call the mirror obtained in this way the $self \ mirror$ of the Ooguri-Vafa space. Clearly this is a fiberwise T duality for nonsingular fibers. However there is no analog of the dual (fiber) metric construction of semiflat mirror symmetry.

Of course taking the self mirror is not  the only way to maintain the fiberwise T duality which is a rather weak requirement anyway. Even if we insist that the mirror must be an Ooguri-Vafa space with the semiflat mirror in the semiflat approximation we still can not determine the singular fibers. In fact any choice of $\Omega_{q}$ that preserves the constraint (87) would work. If we modify the self mirror by only changing $\Omega_{q}$ while keeping (87), the obtained Ooguri-Vafa space is called a $modified \ self \ mirror$.

\section{Ooguri-Vafa Spaces as Local Models of Hitchin's Moduli Spaces}
\subsection{Local Models of Hyperkahler Metrics}
The construction of Fock-Goncharov coordinates as cross-ratios depends on $R$. A key observation due to Gaiotto, Moore and Neitzke is the following proposition.
\begin{thm}\citep{GMN2} For the moduli space $\mathcal{M}$ in section 2.3, the large $R$ asymptotic behavior of $\mathcal{Y}_{\gamma}^{\vartheta}$ is\begin{equation}\mathcal{Y}_{\gamma}^{\vartheta}=\mathcal{Y}_{\gamma}^{sf}(1+O(e^{-cR}))\end{equation}where $c$ is a constant.\end{thm}

The above theorem means that when $R$ goes to $\infty$ the leading order approximation of the hyperkahler metric is given by the semiflat metric. The rescaled semiflat metric $R^{-1}g^{sf}$ collapses to a so-called special Kahler metric on the base of the Hitchin's fibration outside the discriminant locus.

On the other hand the family version of SYZ mirror symmetry conjecture predicts that the SYZ mirror symmetry should be valid for a large complex degeneration of Calabi-Yau varieties and the metric degeneration should have the above behaviors. This fact motivates us to view the large $R$ limit as the differential geometric analog of the so-called large complex degeneration in the standard mirror symmetry (see \citep{L} for discussion). That is why it is quite natural to consider the large $R$ limits of Hitchin's moduli spaces when one studies the SYZ mirror symmetry.

In this paper we are interested in the large $R$ behavior near singular fibers and we need a better approximation than the semiflat one as the latter only takes care of the region outside the discriminant locus. Since we want to use the Ooguri-Vafa space as our local model in this section we only consider the 4-dimensional Hitchin's moduli space $\mathcal{M}^{0}$ and its mirror $\bar{\mathcal{M}}^{0}$. Recall that we have three singular fibers in the Hitchin's fibration of $\mathcal{M}^{0}$. Let $u_{0}$ be one of the three points over which we have singular fibers. Let $\gamma_{e},\gamma_{m}$ be a symplectic basis of the gauge charge lattice. Let $\gamma$ be the vanishing cycle for the degeneration of nearby nonsingular elliptic curves to that singular fiber. Without loss of generality we assume $\gamma=q_{\gamma}\gamma_{e}$. Note that we have the identification $a_{\gamma}=Z_{\gamma}$.
\begin{thm}Suppose conjecture 2.2 is true. As $R\rightarrow\infty$, in a small enough tubular neighborhood of the singular fiber over $u_{0}$ we have \begin{equation}\Omega(\xi)=\Omega^{OV}+O(c_{1}Re^{-c_{2}R})\end{equation} where $c_{1},c_{2}$ depend on the choice of the neighborhood and $c_{2}>0$. Here $\Omega(\xi)$ is the holomorphic symplectic form of the hyperkahler metric form the construction of section 2.3 and $\Omega^{OV}(\xi)$ is the holomorphic symplectic form of the  Ooguri-Vafa space with \begin{equation}\Omega_{2}=\Omega_{-2}=1,\ \ \Omega_{q}=0\ q\neq2 \end{equation}
\end{thm}

\noindent{\bf Proof} First we need the following lemma.
\begin{lem}\citep{GMN3} Suppose conjecture 2.2 is true. Then the holomorphic symplectic form $\Omega(\xi)$ is given by\begin{equation}\Omega(\xi)=\Omega^{sf}(\xi)+\sum_{\gamma^{'}\in\Gamma_{gau}}\Omega_{\gamma^{'}}^{inst}(\xi)+\cdots\end{equation}The ellipses represent a remainder that goes to 0 faster than $\sum_{\gamma^{'}}Re^{-2\pi R|Z_{\gamma^{'}}|}$ (for fixed $Z_{\gamma^{'}}$) as $R\rightarrow \infty$. $\Omega^{sf}$ is the holomorphic symplectic form of the semiflat metric.\begin{equation}\Omega_{\gamma^{'}}^{inst}(\xi)=-\Omega(\gamma^{'}){1\over4\pi^{2} R}{d\mathcal{Y}_{\gamma^{'}}^{sf}(\xi)\over\mathcal{Y}_{\gamma^{'}}^{sf}(\xi)}\wedge[\pi iA^{inst}_{\gamma^{'}}+{1\over2}\pi iV^{inst}_{\gamma^{'}}(\xi^{-1}da_{\gamma^{'}}-\xi d\bar{a}_{\gamma^{'}})]\end{equation}\begin{equation}V^{inst}_{\gamma^{'}}={Rq^{2}_{\gamma^{'}}\over2\pi}\sum_{n>0}e^{in\theta_{\gamma^{'}}}K_{0}(2\pi R|na_{\gamma^{'}}|)\end{equation}Here we have picked a symplectic basis $\{\gamma_{1},\cdots,\gamma_{2r}\}$ ($2r$ is the rank of the gauge charge lattice, for $\mathcal{M}^{0}$ it is $2$) such that $\gamma^{'}=q_{\gamma^{'}}\gamma_{1}$.\begin{equation}A^{inst}_{\gamma^{'}}=-{Rq_{\gamma^{'}}^{2}\over 4\pi}({da_{\gamma^{'}}\over a_{\gamma^{'}}}-{d\bar{a}_{\gamma^{'}}\over \bar{a}_{\gamma^{'}}})\sum_{n>0}e^{in\theta_{\gamma^{'}}}|a_{\gamma^{'}}|K_{1}(2\pi R|na_{\gamma^{'}}|)\end{equation}\end{lem}Once one assumes the conjecture 2.2, we can set up an iteration scheme. The leading term is given by the semiflat one by theorem 6.1. The next to leading term is obtained by using $\mathcal{Y}_{\gamma}^{sf}$ to approximate $\mathcal{Y}_{\gamma}$ on the right hand side of (34). Then a direct calculation in \citep{GMN3} proves the lemma.

Now the terms with $inst$ are exponentially suppressed by $R|a_{\gamma^{'}}|$. Near $u_{0}$ the local coordinate is given by the period $a_{\gamma}(u)$ and $u_{0}$ is  where $a_{\gamma}(u)=0$ due to the vanishing cycle. In a small enough  neighborhood of $u_{0}$ the fiber over $u_{0}$ is the only  singular fiber. In this neighborhood any period $a_{\gamma^{'}}=Z_{\gamma^{'}}$ except the ones with $\gamma^{'}=\pm\gamma$ has a positive lower bound for its norm. Here we do not have to consider $\pm n\gamma$ for $n>1$ because of (33). Therefore in  the sum $\sum_{\gamma^{'}\in\Gamma_{gau}}\Omega_{\gamma^{'}}^{inst}(\xi)$ all but two terms are asymptotically $Re^{-cR}$ for some constant $c$. On the other hand it is easy to check $\Omega^{sf}+\Omega^{inst}_{\gamma}+\Omega^{inst}_{-\gamma}$ is exactly\footnote{The contribution of $-\gamma$ to $A^{inst}$ introduces an extra minus sign.} $\Omega^{OV}(\xi)$ for an Ooguri-Vafa space. In the proof of the next theorem we will show that  for $\mathcal{M}^{0}$ \begin{equation}\Omega_{q_{\gamma}}=\Omega_{-q_{\gamma}}=1,\ \ \Omega_{q}=0\ otherwise\end{equation}For $\mathcal{M}^{0}$ we have $q_{\gamma}=2$ as two circles (each of them representing $\gamma_{e}$) collapse to two points in the singular fiber over $u_{0}$.\\

Since the hyperkahler structure is determined by $\Omega(\xi)$, the theorem above tells us that in a small neighborhood of $u_{0}$ the Ooguri-Vafa metric is a good approximation of the hyperkahler structure of $\mathcal{M}^{0}$ in the large $R$ limit (large complex limit). Note that by the general discussion of Ooguri-Vafa spaces in section 4.1 there are two singularities corresponding to the two singularities in the singular fibers.\footnote{In \citep{FW}, it is pointed out that the two singularities are not singularities of the total spaces while for the Ooguri-Vafa space the two singularities are singular in the total space. There is no contradiction as the Ooguri-Vafa space is an approximation.}

The description of the Ooguri-Vafa space in section 4.2 in terms of twistor coordinates exhibits nontrivial discontinuous jumps and trivial 4d wall crossing. We want to see how this matches the situation of the full Hitchin's moduli space $\mathcal{M}^{0}$.
\begin{thm}The discontinuous jumps of Ooguri-Vafa spaces can be identified with discontinuous jumps in the Gaiotto-Moore-Neitzke constructions of section 2.3.\end{thm}

\noindent{\bf Proof}  Our Riemann surface $C$ is a nonsingular elliptic curve and the spectral curves are also elliptic curves.  The equation of the spectral curve is $\lambda^{2}=-det\ \varphi=Tr\ \varphi^{2}$. The explicit form of $Tr\ \varphi^{2}$ is determined in \citep{FW} from which it is easy to see that  there are two simple zeroes (denoted by $z_{1}, z_{2}$) of the quadratic differential $\lambda^{2}$. Each spectral curve is a double cover of $C$ with two branch points. There is one regular singularity (for Hitchin's equations) denoted by $z_{0}$. According to the description of the WKB triangulations in section 2.3 there are three separating trajectories for each simple zero and they all have $z_{0}$ as one end. To build a triangulation we need to pick one generic trajectory from each component of the foliation separated by these separating trajectories.

Let us consider the three separating trajectories associated to $z_{1}$. The union of any two of them must be a topologically nontrivial circle on the torus $C$ because otherwise the union bounds an disk area on $C$ in which the trajectories all start and end at $z_{0}$. As we go to $z_{0}$ transversely to the trajectories we encounter a  behavior of generic trajectories contradicting  the local classification in section 2.3. As a result, when we pick three generic trajectories to form a triangle containing $z_{1}$ each one of them is a topologically nontrivial circle. Denote these three circles by $a,b,c$. Choose a symplectic basis $\{\gamma_{1},\gamma_{2}\}$ of $C$.  Suppose one of them (say $a$) represents $\gamma_{1}$. This implies that the sum of the two relevant separating trajectories represents $\gamma_{1}$. The other two sides of the triangle $abc$ must be topologically nontrivial circles, too. So far we have located three edges of a WKB triangulation. By formula (25) that is all we should have. That means that $a,b,c$ are also representatives of generic trajectories separated by three separating trajectories emanating from the other simple zero $z_{2}$. In particular  the sum of two of three separating trajectories emanating from $z_{2}$ that correspond to $a$ also represents $\gamma_{1}$. The above construction seems to depend on an additional assumption that $a$ represents a topologically nontrivial circle winding the torus only once (i.e. $\gamma_{1}$). Note that we only need the existence of such configuration for some $\lambda^{2}$ and $\vartheta$. So we can choose a $\mathbf{Z}_{2}$ (after all we can realize $C$ as a double cover of $CP^{1}$ with four branch points) symmetric foliation such that $z_{0}$ is fixed while $z_{1},z_{2}$ are permuted by the $\mathbf{Z}_{2}$ action and $\gamma_{1}$ is $\mathbf{Z}_{2}$ invariant. Then there is a $\mathbf{Z}_{2}$ invariant generic trajectory representing $\gamma_{1}$.

Therefore we have a WKB triangulation with one vertex, three edges and two triangles. Triangle 1 contains $z_{1}$ while triangle 2 contains $z_{2}$. There are three charges corresponding to the three edges $a,b,c$ by the construction in section 2.3. They generate the charge lattice $\Gamma$ and in particular the gauge charge lattice $\Gamma_{gau}$. If we take the edge $a$ then the above construction produces an cycle $\gamma_{a}$ in the spectral curve which is twice of a topologically nontrivial circle winding around the spectral curve (a torus) once. Now recall that the vanishing cycle at $u_{0}$ is $2\gamma_{e}$. So without loss of generality we can assume that the choice of $\gamma_{1}$ in $H_{1}(C)$ is such that \begin{equation}\gamma_{a}=2\gamma_{e}\end{equation}The BPS number $\Omega(\gamma_{a})$ is $1$ as this is a flip (it corresponds to a finite trajectory connecting two different simple zeroes). So we recover (138) after using (111). The form of $Z_{\gamma_{e}}$ in section 4.1 is by definition while the form of $Z_{\gamma_{m}}$ is from the standard result of Picard-Fuchs equations of periods which gives logarithmic divergence and Picard-Lefschetz transformations (which fixes $\Delta$ to be 4). Of course this result also follows from the Gibbons-Hawking ansatz itself once one identifies the precise Ooguri-Vafa metric with (138). Note that here when we relate $Z_{\gamma_{m}}$ in the Ooguri-Vafa space with the central charge (period) one only keeps the leading order approximation. That is why in section 4.3 we are not very interested in the higher order contributions of $W$ (i.e. $W^{analytic}$). It may not be consistent to extract higher order information in the Ooguri-Vafa approximation.

Since we have identified $\Omega$ in the Ooguri-Vafa space with the BPS numbers in the Gaiotto-Moore-Neitzke construction the identification of discontinuities follows. \\

What happens to the 4d wall crossing formula? To completely determine all 4d wall crossing formulas for $\mathcal{M}^{0}$ is an intricate problem but is hardly relevant in our local model discussion as the 4d wall crossing formula in the Ooguri-Vafa space is trivial. In fact the moral here is that the local model simplifies the behavior of the hyperkahler metric considerably. For the 4-dimensional $\mathcal{M}^{0}$ Ooguri-Vafa spaces take care of all possible local models if we keep only the leading order approximation. It is likely that for higher dimensional Hitchin's moduli spaces more information can be preserved in the local models that are more complicated than (the straightforward generalization of) Ooguri-Vafa models and some truncations of 4d wall crossing formulas will emerge from these local models.

\subsection{Mirror Symmetry and Local Models}

Recall that the SYZ mirror of $\mathcal{M}^{0}$ (for some fixed large enough $R$) is $\bar{\mathcal{M}}^{0}=\mathcal{M}^{0}/Q$. So a hyperkahler metric  of $\bar{\mathcal{M}}^{0}$ is induced from the one on $\mathcal{M}^{0}$. The analog of conjecture 2.1 is\begin{con}For some fixed large enough $R$, $\bar{\mathcal{M}}^{0}$ endowed with the hyperkahler quotient metric is isometric to $\mathcal{M}^{0}/Q$ endowed with the quotient metric of the Gaiotto-Moore-Neitzke hyperkahler metric.\end{con}

Even if the conjecture is true it is still very hard to discuss the metrical aspect of the SYZ mirror symmetry for two reasons. First the Gaiotto-Moore-Neitzke construction is not explicit. Second there is no direct construction of hyperkahler metrics of $SO(3)$ Hitchin's moduli spaces analogous to the  Gaiotto-Moore-Neitzke construction.

Now let us consider the Ooguri-Vafa local model. As pointed out in section 5.3  there is no canonical choice of an Ooguri-Vafa type SYZ mirror. However here we have one additional requirement. The dual fiber over $u_{0}$ has one singularity. So only for $q=1$ there is nonzero $\Omega_{q}$ in the sum. What is the value of $\Omega_{1}$? The lack of a geometric construction of BPS numbers gives us no hints. However the monodromy around the singular fiber suggests that $\Omega_{1}=2$ so that $\Delta=2$ (suggested by conjecture 6.1).

\begin{con}Let $\Omega(\xi)$ be the holomorphic symplectic form of the hyperkahler metric of $\bar{\mathcal{M}}^{0}$ (the one isometric to $\mathcal{M}^{0}/Q$).  As $R\rightarrow\infty$, in a small enough tubular neighborhood of the singular fiber over $u_{0}$ we have \begin{equation}\Omega(\xi)=\Omega^{OV}+O(c_{1}Re^{-c_{2}R})\end{equation} where $c_{1},c_{2}$ depend on the choice of the neighborhood and $c_{2}>0$. Here $\Omega^{OV}(\xi)$ is the holomorphic symplectic form of the  Ooguri-Vafa space with \begin{equation}\Omega_{1}=\Omega_{-1}=2,\ \ \Omega_{q}=0\ q\neq\pm1 \end{equation}\end{con}

Now that we have some understandings of the local model of the hyperkahler structure we may want to study local models of hyperholomorphic connections and SYZ mirror symmetry of them. We will discuss this problem in another paper.

\end{document}